\documentclass[12pt]{article}          \usepackage{amsfonts} \usepackage{amssymb} \usepackage{latexsym} \usepackage{graphicx} \usepackage{psfrag} \usepackage{amsbsy}        
\newtheorem{thm}{Theorem}

\newcounter{claimcount}[thm]  
\newtheorem{prop}[thm]{Proposition} 

\newtheorem{lemma}[thm]{Lemma} 

\newtheorem{cor}[thm]{Corollary}

\newcommand{\bprf}[1][Proof:]{\begin{list}{} 			{\setlength{\leftmargin}{1em} 			\setlength{\rightmargin}{0em}}                         \item {\bf \hspace{-1em}  #1 \ \ }} 
\newcommand{\entr}[2][\mu]{h_{#1}\left(#2\right)}
\newcommand{\Entr}[3][\mu]{H_{#1}\left(#2\left|#3\right.\right)}
\newcommand{\Expct}[3][\mu]{{\mathbb{E}}_{#1}\left[#2\left|#3\right.\right]}

\begin{document}

\title{Invariant measures for bipermutative cellular automata}

\author{
Marcus Pivato\thanks{Partially supported by NSERC Canada.}\\ {\em Department of Mathematics, Trent University} \\ {\tt pivato@xaravve.trentu.ca}}

\maketitle

\begin{abstract} A {\em right-sided, nearest neighbour cellular
automaton} (RNNCA) is a continuous transformation
$\Phi:{\mathcal{ A}}^{\mathbb{Z}}{{\longrightarrow}}{\mathcal{ A}}^{\mathbb{Z}}$ determined by a {\em local rule}
$\phi:{\mathcal{ A}}^{\{0,1\}}{{\longrightarrow}}{\mathcal{ A}}$ so that, for any ${\mathbf{ a}}\in{\mathcal{ A}}^{\mathbb{Z}}$ and any
$z\in{\mathbb{Z}}$,\quad $\Phi({\mathbf{ a}})_z \ = \ \phi(a_{z},a_{z+1})$.  We say that
$\Phi$ is {\em bipermutative} if, for any choice of $a\in{\mathcal{ A}}$, the map
${\mathcal{ A}}\ni b \mapsto \phi(a,b)\in{\mathcal{ A}}$ is bijective, and also, for any
choice of $b\in{\mathcal{ A}}$, the map ${\mathcal{ A}}\ni a \mapsto \phi(a,b)\in{\mathcal{ A}}$ is
bijective.

  We characterize the invariant measures of bipermutative RNNCA.
First we introduce the equivalent notion of a {\em quasigroup CA}, to
expedite the construction of examples.  Then we characterize
$\Phi$-invariant measures when ${\mathcal{ A}}$ is a (nonabelian) group, and
$\phi(a,b) = a\cdot b$.  Then we show that, if $\Phi$ is any
bipermutative RNNCA, and $\mu$ is $\Phi$-invariant, then $\Phi$
must be $\mu$-almost everywhere $K$-to-1, for some constant $K$.
We use this to characterize invariant measures when
${\mathcal{ A}}^{\mathbb{Z}}$ is a group shift and $\Phi$ is an endomorphic CA.

\end{abstract}
\section{Introduction}

 If ${\mathcal{ A}}$ is a (discretely topologized) finite set, then ${\mathcal{ A}}^{\mathbb{Z}}$
is compact in the Tychonoff topology.  Let
$ {{{\boldsymbol{\sigma}}}^{}} :{\mathcal{ A}}^{\mathbb{Z}}{{\longrightarrow}}{\mathcal{ A}}^{\mathbb{Z}}$ be the shift map: $ {{{\boldsymbol{\sigma}}}^{}} ({\mathbf{ a}}) \ = \
{\left[b_z  |_{z\in{\mathbb{Z}}}^{} \right]}$, where $b_z = a_{z-1}$, \ $\forall
z\in{\mathbb{Z}}$.  A {\bf cellular automaton} (CA) is a continuous map
$\Phi:{\mathcal{ A}}^{\mathbb{Z}}{{\longrightarrow}}{\mathcal{ A}}^{\mathbb{Z}}$ which commutes with $ {{{\boldsymbol{\sigma}}}^{}} $.
Equivalently, $\Phi$ is defined by a {\bf local rule}
$\phi:{\mathcal{ A}}^{\left[ -\ell...r \right]}{{\longrightarrow}}{\mathcal{ A}}$ (for some $\ell, r \geq 0$)
so that, for any ${\mathbf{ a}}\in{\mathcal{ A}}^{\mathbb{Z}}$ and any $z\in{\mathbb{Z}}$,\quad
$\Phi({\mathbf{ a}})_z \ = \ \phi(a_{z-\ell},\ldots,a_{z+r})$.
  We say $\Phi$ is {\bf right-permutative} if, for any fixed
${\mathbf{ a}}\in{\mathcal{ A}}^{\left[ -\ell...r \right)}$, the map ${\mathcal{ A}}\ni b
\mapsto \phi({\mathbf{ a}},b)\in {\mathcal{ A}}$ is bijective.  Likewise, $\Phi$ is {\bf
left-permutative} if, for any fixed ${\mathbf{ b}}\in{\mathcal{ A}}^{\left( -\ell...r \right]}$, the map ${\mathcal{ A}}\ni a \mapsto \phi(a,{\mathbf{ b}})\in {\mathcal{ A}}$ is bijective, and $\Phi$ is {\bf
bipermutative} if it is both left- and right-permutative.

  \medskip         \refstepcounter{thm} {\bf Example \thethm:}  \setcounter{enumi}{\thethm} \begin{list}{(\alph{enumii})}{\usecounter{enumii}} 			{\setlength{\leftmargin}{0em} 			\setlength{\rightmargin}{0em}}   
  \item \label{X:biperm.CA.a} 
If $({\mathcal{ A}},+)$ is an abelian group, $\ell=0$ and $r=1$, and $\phi(a_0,a_1) =
a_0 + a_1$, then $\Phi$ is a called a {\em nearest neighbour addition} CA,
and is bipermutative.
 
  \item  \label{X:biperm.CA.b}
  If ${\mathcal{ A}}={{\mathbb{Z}}_{/p}}$, and let $c_0,c_1\in{\left[ 1..p \right)}$
be constants.  If $\phi(a_0,a_1) = c_0a_0 + c_1a_1$,
then $\Phi$ is a called a {\em Ledrappier} CA,
and is bipermutative.
 	\hrulefill\end{list}   			
  We say that $\Phi$ is a {\bf right-sided, nearest neighbour}
cellular automaton (RNNCA) if  $\ell=0$ and $r=1$
(as in Examples \ref{X:biperm.CA.a} and
\ref{X:biperm.CA.b}).  It is easy to show:

\begin{lemma}{\sf }   Let $\Phi:{\mathcal{ A}}^{\mathbb{Z}}{{\longrightarrow}}{\mathcal{ A}}^{\mathbb{Z}}$ be a CA and let
${\mathcal{ B}}={\mathcal{ A}}^{\ell+r}$.  There is an RNNCA $\Gamma:{\mathcal{ B}}^{\mathbb{Z}}{{\longrightarrow}}{\mathcal{ B}}^{\mathbb{Z}}$
so that the the topological dynamical system $({\mathcal{ A}}^{\mathbb{Z}},\Phi)$ is
isomorphic to the system $({\mathcal{ B}}^{\mathbb{Z}},\Gamma)$.

  Furthermore $\left( \ \rule[-0.5em]{0em}{1em}       \begin{minipage}{40em}       \begin{tabbing}         $\Phi$ is bipermutative        \end{tabbing}      \end{minipage} \ \right)
\iff \left( \ \rule[-0.5em]{0em}{1em}       \begin{minipage}{40em}       \begin{tabbing}         $\Gamma$ is bipermutative        \end{tabbing}      \end{minipage} \ \right)$.
\hrulefill\ensuremath{\Box}  \end{lemma}
  Let $\lambda$ be the uniform Bernoulli measure on ${\mathcal{ A}}^{\mathbb{Z}}$.  Thus,
for any $c_1,\ldots,c_M\in{\mathcal{ A}}$, and any $z_1,\ldots,z_M\in{\mathbb{Z}}$,
\[
  \lambda{\left\{ {\mathbf{ a}}\in{\mathcal{ A}}^{\mathbb{Z}} \; ; \; a_{z_1}=c_1,\ldots,a_{z_M}=c_M \right\} } 
\quad=\quad \frac{1}{|{\mathcal{ A}}|^M}.
\]
  Any permutative CA is surjective, and any surjective CA preserves
$\lambda$ \cite{HedlundCA}.  What other $ {{{\boldsymbol{\sigma}}}^{}} $-invariant measures on
${\mathcal{ A}}^{\mathbb{Z}}$ are also invariant under permutative CA?  Let $\entr{\Phi}$
denote the {\bf entropy} \cite[\S5.2]{Petersen} of the measure-preserving
dynamical system $({\mathcal{ A}}^{\mathbb{Z}},\Phi,\mu)$.  Host, Maass, and Martinez
\cite{HostMaassMartinez} have shown:

\begin{prop}{\sf \label{Host.Maass.Martinez}}  
 Let ${\mathcal{ A}}={{\mathbb{Z}}_{/p}}$, where $p$ is prime, and
let $\Phi:{{\mathcal{ A}}^{\mathbb{Z}}}{{\longrightarrow}}{{\mathcal{ A}}^{\mathbb{Z}}}$ be a Ledrappier CA.  
Suppose $\mu$ is a measure which is $\Phi$-invariant and $ {{{\boldsymbol{\sigma}}}^{}} $-ergodic.
If $\entr{\Phi}>0$, then $\mu=\lambda$.\hrulefill\ensuremath{\Box}
 \end{prop}

  This paper provides generalizations of Proposition \ref{Host.Maass.Martinez}
to a variety of contexts.
  In \S\ref{S:loop}, we introduce quasigroups, which provide a convenient
formulation of bipermutative RNNCA as {\em quasigroup cellular automata},
and suggest a natural generalization of Proposition
\ref{Host.Maass.Martinez} (Conjecture \ref{subloop.conjecture}).
In \S\ref{S:group} we characterize invariant measures for {\em
nearest-neighbour multiplication} CA (when ${\mathcal{ A}}$ is a nonabelian
group), and construct an explicit counterexample to Conjecture
\ref{subloop.conjecture}. In \S\ref{S:degree} we will extend the
method of \cite{HostMaassMartinez} to prove that, if $\mu$ is $\Phi$-invariant,
then there is some $K\leq|{\mathcal{ A}}|$ so that $\Phi$ is $K$-to-1 {\rm ($\mu$-\ae)}
\ (Theorem \ref{loop.ca.K.2.1}).  In
\S\ref{S:eca} we will provide a generalization of Proposition
\ref{Host.Maass.Martinez} to {\em endomorphic} CA on group shifts
(Theorem \ref{biperm.eca.measure}).

\paragraph*{Notation:}  If $\mu$ is a measure, then `${{\forall}_{\mu}\;}  x$'
means `for $\mu$-almost all $x$', and `$\mu$-\ae' means `$\mu$-almost
everywhere'.  If ${\mathbf{ U}}$ and ${\mathbf{ V}}$ are measurable sets, then  `${\mathbf{ U}} \ \raisebox{-1ex}{$\stackrel{\displaystyle\subset}{\scriptstyle{\mathrm{\mu}}}$} \ {\mathbf{ V}}$'
means $\mu[{\mathbf{ U}}\setminus{\mathbf{ V}}]=0$, and `${\mathbf{ U}} \
\displaystyle\raisebox{-0.6ex}{$\overline{\overline{{\scriptstyle{\mathrm{\,\mu\,}}}}}$} \ {\mathbf{ V}}$' means ${\mathbf{ U}} \ \raisebox{-1ex}{$\stackrel{\displaystyle\subset}{\scriptstyle{\mathrm{\mu}}}$} \ {\mathbf{ V}}$ and
${\mathbf{ V}} \ \raisebox{-1ex}{$\stackrel{\displaystyle\subset}{\scriptstyle{\mathrm{\mu}}}$} \ {\mathbf{ U}}$.  If ${\mathfrak{ S}}$ is a sigma algebra
and ${\mathbf{ U}}$ is a measurable set, then $\Expct{{\mathbf{ U}}}{{\mathfrak{ S}}}$ is the
conditional expectation of ${\mathbf{ U}}$ given ${\mathfrak{ S}}$.

\section{\label{S:loop}\label{S:dualCA} Quasigroup Cellular Automata}

  A {\bf quasigroup} \cite{Pflugfelder} is a finite set ${\mathcal{ A}}$ equipped
with a binary operation `$*$' which has the left- and
right-cancellation properties.  In other words, for any $a,b,c\in{\mathcal{ A}}$,
\[
 \left( \ \rule[-0.5em]{0em}{1em}       \begin{minipage}{40em}       \begin{tabbing}         $a*b = a*c$        \end{tabbing}      \end{minipage} \ \right) \ensuremath{\Longrightarrow} \left( \ \rule[-0.5em]{0em}{1em}       \begin{minipage}{40em}       \begin{tabbing}         $b=c$        \end{tabbing}      \end{minipage} \ \right),
\quad\mbox{\ and \ }\quad
 \left( \ \rule[-0.5em]{0em}{1em}       \begin{minipage}{40em}       \begin{tabbing}         $b*a = c*a$        \end{tabbing}      \end{minipage} \ \right) \ensuremath{\Longrightarrow} \left( \ \rule[-0.5em]{0em}{1em}       \begin{minipage}{40em}       \begin{tabbing}         $b=c$        \end{tabbing}      \end{minipage} \ \right).
\]
  If we identify ${\mathcal{ A}}$ with ${\left[ 1..N \right]}$ in some arbitrary way, then
the `multiplication table' for $*$ is the $N\times N$ matrix ${\mathbf{ M}}^* =
[m_{i,j}]_{i,j=1}^N$ where $m_{i,j} = i * j$.  We say ${\mathbf{ M}}^*$ is
a {\bf Latin square} \cite{DenesKeedwell} if every column and every row of ${\mathbf{ M}}^*$ contains each element of
${\left[ 1..N \right]}$ exactly once.  It follows:
\[
\left( \ \rule[-0.5em]{0em}{1em}       \begin{minipage}{40em}       \begin{tabbing}         $({\mathcal{ A}},*)$ is a quasigroup        \end{tabbing}      \end{minipage} \ \right) \iff
\left( \ \rule[-0.5em]{0em}{1em}       \begin{minipage}{40em}       \begin{tabbing}         ${\mathbf{ M}}^*$ is a Latin square        \end{tabbing}      \end{minipage} \ \right).
\]
  Note that the operator `$*$' is not necessarily associative.  Indeed,
it is easy to show:
\[
 \left( \ \rule[-0.5em]{0em}{1em}       \begin{minipage}{40em}       \begin{tabbing}         `$*$' is associative        \end{tabbing}      \end{minipage} \ \right) \iff \left( \ \rule[-0.5em]{0em}{1em}       \begin{minipage}{40em}       \begin{tabbing}         $({\mathcal{ A}},*)$ is a group        \end{tabbing}      \end{minipage} \ \right).
\]
  A {\bf quasigroup cellular automaton} (QGCA) is a right-sided, nearest
neighbour cellular automaton $\Phi:{{\mathcal{ A}}^{\mathbb{Z}}}{{\longrightarrow}}{{\mathcal{ A}}^{\mathbb{Z}}}$ with local 
rule $\phi:{\mathcal{ A}}^{\{0,1\}}{{\longrightarrow}}{\mathcal{ A}}$ given: $\phi(a_0,a_1) = a_0 * a_1$,
where `$*$' is a quasigroup operation.  For example, any Ledrappier automaton
is a QGCA.  It follows:

\begin{prop}{\sf }  
 $ \left( \ \rule[-0.5em]{0em}{1em}       \begin{minipage}{40em}       \begin{tabbing}         $\Phi$ is a bipermutative RNNCA        \end{tabbing}      \end{minipage} \ \right)
\quad \iff \quad 
\left( \ \rule[-0.5em]{0em}{1em}       \begin{minipage}{40em}       \begin{tabbing}         $\Phi$ is a quasigroup CA        \end{tabbing}      \end{minipage} \ \right)$.\hrulefill\ensuremath{\Box}
 \end{prop}

The obvious generalization of Proposition \ref{Host.Maass.Martinez}
fails for arbitrary quasigroup CA.  If ${\mathcal{ B}}\subset{\mathcal{ A}}$, then we call ${\mathcal{ B}}$
a {\bf subquasigroup}  (and write `${\mathcal{ B}}\prec{\mathcal{ A}}$') if ${\mathcal{ B}}$ is closed under
the `$*$' operation.  

\begin{lemma}{\sf }   If $\Phi:{{\mathcal{ A}}^{\mathbb{Z}}}{{\longrightarrow}}{{\mathcal{ A}}^{\mathbb{Z}}}$ is a QGCA, and ${\mathcal{ B}}\prec{\mathcal{ A}}$,
then ${\mathcal{ B}}^{\mathbb{Z}}$ is a $\Phi$-invariant subshift.
If $\mu$ is the uniform Bernoulli measure on ${\mathcal{ B}}^{\mathbb{Z}}$, then
$\mu$ is $\Phi$-invariant and $ {{{\boldsymbol{\sigma}}}^{}} $-ergodic.  If $|{\mathcal{ B}}|=K$, then
$\entr{\Phi} \ = \ \log(K)$ and $\Phi$ is $K$-to-1 {\rm ($\mu$-\ae)}.\hrulefill\ensuremath{\Box}
 \end{lemma}

  If ${\mathcal{ B}}\prec{\mathcal{ A}}$ and $({\mathcal{ A}},*)$ is a finite group, then ${\mathcal{ B}}$ is a
subgroup.  Thus, if $|{\mathcal{ A}}|$ is prime, then ${\mathcal{ A}}$ can't have nontrivial
subquasigroups.  However, other prime cardinality quasigroups can:

        \refstepcounter{thm}                     \begin{list}{} 			{\setlength{\leftmargin}{1em} 			\setlength{\rightmargin}{0em}}        \item {\bf Example \thethm:} Let ${\mathcal{ D}} = \{a_1,a_2; \ b_1,b_2; \ c_1,c_2,c_3\}$;
thus, $|{\mathcal{ D}}|=7$ is prime.  Let $*$ have the following multiplication
table:
\[
\begin{array}{c||c|c|c|c|c|c|c|}
*   & a_1\hspace{0.6em}  & a_2 & b_1 \hspace{0.6em} & b_2 & c_1 & c_2 & c_3 
\\ \hline  \hline \cline{2-3}\cline{2-3}
a_1 & \multicolumn{2}{||c||}{\begin{array}{c|c} a_1 & a_2\end{array}} & c_1 & c_2 & b_2 & b_1 & c_3 \\ \hline
a_2 & \multicolumn{2}{||c||}{\begin{array}{c|c} a_2 & a_1 \end{array}} & c_2 & c_1 & b_1 & c_3 & b_2 \\ \hline \cline{2-5} \cline{2-5}
b_1 & c_1 & c_3 & \multicolumn{2}{||c||}{\begin{array}{c|c}b_1 & b_2 \end{array}}& c_2 & a_1 & a_2 \\ \hline
b_2 & c_3 & c_1 & \multicolumn{2}{||c||}{\begin{array}{c|c}b_2 & b_1 \end{array}} & a_1 & a_2 & c_2 \\ \hline \cline{4-5} \cline{4-5}
c_1 & b_1 & b_2 & c_3 & a_1 & a_2 & c_2 & c_1 \\ \hline
c_2 & b_2 & c_2 & a_1 & a_2 & c_3 & c_1 & b_1 \\ \hline
c_3 & c_2 & b_1 & a_2 & c_3 & c_1 & b_2 & a_1 \\ \hline
\end{array}
\]
  Clearly, the quasigroup $({\mathcal{ D}},*)$ has two subquasigroups:
${\mathcal{ A}}=\{a_1,a_2\}$ and ${\mathcal{ B}}=\{b_1,b_2\}$.
                   \hrulefill  \end{list}   			

  This suggests that the correct generalization of Proposition \ref{Host.Maass.Martinez} is:

\paragraph*{\refstepcounter{thm}\label{subloop.conjecture} Conjecture \thethm:}
{\em
 Let $\Phi:{{\mathcal{ A}}^{\mathbb{Z}}}{{\longrightarrow}}{{\mathcal{ A}}^{\mathbb{Z}}}$ be a QGCA.
If $\mu$ is a $\Phi$-invariant and $ {{{\boldsymbol{\sigma}}}^{}} $-ergodic measure,
and $\entr{\Phi}>0$, then $\mu$ is the uniform measure on ${\mathcal{ B}}^{\mathbb{Z}}$,
for some ${\mathcal{ B}}\prec{\mathcal{ A}}$. }

\medskip

  Conjecture \ref{subloop.conjecture} is false, as we will show with
Example (\ref{counterexample}) of \S\ref{S:group}.

\paragraph*{Unilateral vs. Bilateral Cellular Automata:}

 Any right-sided CA $\Phi:{\mathcal{ A}}^{\mathbb{Z}}{{\longrightarrow}}{\mathcal{ A}}^{\mathbb{Z}}$ induces a {\bf
unilateral} CA $\widetilde\Phi:{\mathcal{ A}}^{\mathbb{N}}{{\longrightarrow}}{\mathcal{ A}}^{\mathbb{N}}$ with the same local
rule.  Any $\Phi$-invariant measure on ${\mathcal{ A}}^{\mathbb{Z}}$ projects to a
$\widetilde\Phi$-invariant measure on ${\mathcal{ A}}^{\mathbb{N}}$; conversely, any
$\widetilde\Phi$- and $ {{{\boldsymbol{\sigma}}}^{}} $-invariant measure on ${\mathcal{ A}}^{\mathbb{N}}$ extends to
a unique $(\Phi, {{{\boldsymbol{\sigma}}}^{}} )$-invariant measure on ${\mathcal{ A}}^{\mathbb{Z}}$.
In what follows, we will abuse notation and write $\widetilde\Phi$ as $\Phi$.
Thus, Conjecture \ref{subloop.conjecture} is equivalent to:

\paragraph*{Conjecture $\widetilde{\ref{subloop.conjecture}}$:}
{\em
 Let $\Phi:{\mathcal{ A}}^{\mathbb{N}}{{\longrightarrow}}{\mathcal{ A}}^{\mathbb{N}}$ be a (unilateral) QGCA.
If $\mu$ is $\Phi$-invariant and $ {{{\boldsymbol{\sigma}}}^{}} $-ergodic,
and $\entr{\Phi}>0$, then $\mu$ is the uniform measure on ${\mathcal{ B}}^{\mathbb{N}}$,
for some ${\mathcal{ B}}\prec{\mathcal{ A}}$. }

\paragraph*{Dual Cellular Automata:} There is a well-known
conjugacy between any right-permutative unilateral CA and a full
shift.  Define $\Xi:{\mathcal{ A}}^{\mathbb{N}} {{\longrightarrow}} {\mathcal{ A}}^{\mathbb{N}}$ by \ $ \Xi({\mathbf{ a}}) \ = \
\left[a_0, \ \Phi({\mathbf{ a}})_0, \ \Phi^2({\mathbf{ a}})_0, \ \Phi^3({\mathbf{ a}})_0, \ \ldots
\right]$.

\begin{lemma}{\sf \label{phi.to.shift.conj}}  
  If $\Phi$ is right-permutative, then 
 $\Xi$ is a topological conjugacy from the dynamical system
$({\mathcal{ A}}^{\mathbb{N}},\Phi)$ to the system $({\mathcal{ A}}^{\mathbb{N}}, {{{\boldsymbol{\sigma}}}^{}} )$
(ie. $\Xi$ is a homeomorphism and $\Xi \circ \Phi \ = \  {{{\boldsymbol{\sigma}}}^{}} \circ\Xi$).
\hrulefill\ensuremath{\Box}
 \end{lemma}

 Let $({\mathcal{ A}},*)$ be a quasigroup.  The {\bf dual} quasigroup is the set
${\mathcal{ A}}$ equipped with binary operator $\,\widehat{*}\,$ defined: $a \,\widehat{*}\,
b \ = \ c$, where $c$ is the unique element in ${\mathcal{ A}}$ such that $a*c =
b$.  If $({\mathcal{ A}},*)$ is a group, then $a\,\widehat{*}\, b = a^{-1}* b$.  If
$\Phi:{\mathcal{ A}}^{\mathbb{Z}}{{\longrightarrow}}{\mathcal{ A}}^{\mathbb{Z}}$ is a QGCA (with local map
$\phi(a,b)=a*b$), then the {\bf dual} of $\Phi$ is the CA
${\widehat{\Phi }}:{\mathcal{ A}}^{\mathbb{Z}}{{\longrightarrow}}{\mathcal{ A}}^{\mathbb{Z}}$ having local map ${\widehat{\phi }}(a,b)
= a\,\widehat{*}\, b$.

\begin{lemma}{\sf \label{shift.to.xi.conj}}  
Let  $({\mathcal{ A}},*)$ be a quasigroup and let $\Phi$ be the corresponding QGCA.  Then:
 \setcounter{enumi}{\thethm} \begin{list}{{\bf (\alph{enumii})}}{\usecounter{enumii}} 			{\setlength{\leftmargin}{0em} 			\setlength{\rightmargin}{0em}}
  \item $({\mathcal{ A}},\widehat*)$ is a quasigroup, and
${\widehat{\Phi }}$ is a QGCA. The dual of $\,\widehat{*}\,$ is $*$; \ the dual of ${\widehat{\Phi }}$ is $\Phi$.

  \item $\Xi$ is a topological conjugacy from the dynamical system
$({\mathcal{ A}}^{\mathbb{N}}, {{{\boldsymbol{\sigma}}}^{}} )$ to the system $({\mathcal{ A}}^{\mathbb{N}},{\widehat{\Phi }})$,
so that we have the following commuting cube:

\begin{center}
\psfrag{A}[][]{${\mathcal{ A}}^{\mathbb{N}}$}
\psfrag{f}[][]{$\Phi$}
\psfrag{s}[][]{$ {{{\boldsymbol{\sigma}}}^{}} $}
\psfrag{p}[][]{$\widehat{\Phi}$}
\psfrag{X}[][]{$\Xi$}
\includegraphics[scale=0.7]{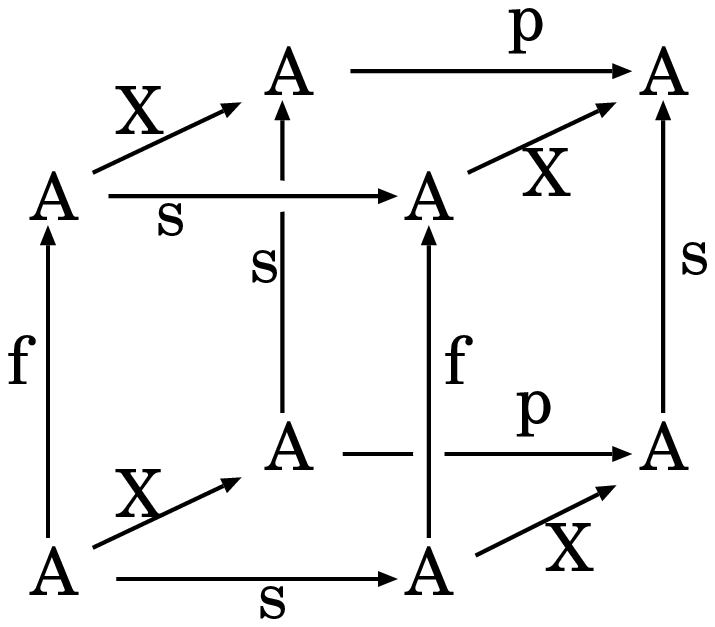}
\end{center}

  \item If ${\mathcal{ B}}\subset{\mathcal{ A}}$, then $\left( \ \rule[-0.5em]{0em}{1em}       \begin{minipage}{40em}       \begin{tabbing}         $({\mathcal{ B}},*)\prec({\mathcal{ A}},*)$        \end{tabbing}      \end{minipage} \ \right)\iff\left( \ \rule[-0.5em]{0em}{1em}       \begin{minipage}{40em}       \begin{tabbing}         $({\mathcal{ B}},\,\widehat{*}\,) \prec ({\mathcal{ A}},\,\widehat{*}\,)$        \end{tabbing}      \end{minipage} \ \right)$.
\end{list}
  Let $\mu$ be a measure on ${\mathcal{ A}}^{\mathbb{Z}}$, and let ${\widehat{\mu}} = \Xi(\mu)$.  Then:
 \setcounter{enumi}{\thethm} \begin{list}{{\bf (\alph{enumii})}}{\usecounter{enumii}} 			{\setlength{\leftmargin}{0em} 			\setlength{\rightmargin}{0em}}
  \item $\left( \ \rule[-0.5em]{0em}{1em}       \begin{minipage}{40em}       \begin{tabbing}         $\mu$ is $\Phi$-invariant        \end{tabbing}      \end{minipage} \ \right)\iff
\left( \ \rule[-0.5em]{0em}{1em}       \begin{minipage}{40em}       \begin{tabbing}         ${\widehat{\mu}}$ is $ {{{\boldsymbol{\sigma}}}^{}} $-invariant        \end{tabbing}      \end{minipage} \ \right)$.

  \item $\left( \ \rule[-0.5em]{0em}{1em}       \begin{minipage}{40em}       \begin{tabbing}         $\mu$ is $ {{{\boldsymbol{\sigma}}}^{}} $-ergodic        \end{tabbing}      \end{minipage} \ \right)\iff
\left( \ \rule[-0.5em]{0em}{1em}       \begin{minipage}{40em}       \begin{tabbing}         ${\widehat{\mu}}$ is ${\widehat{\Phi }}$-ergodic        \end{tabbing}      \end{minipage} \ \right)$.

  \item If $\mu$ is $\Phi$- and $ {{{\boldsymbol{\sigma}}}^{}} $- invariant, then
 $h(\Phi,\mu) \ = \ h({\widehat{\Phi }},{\widehat{\mu}})
\ = \ h( {{{\boldsymbol{\sigma}}}^{}} ,\mu) \ = \ h( {{{\boldsymbol{\sigma}}}^{}} ,{\widehat{\mu}})$.\hrulefill\ensuremath{\Box}
\end{list}
 \end{lemma}

Thus, Conjecture $\widetilde{\ref{subloop.conjecture}}$ is equivalent to:

\paragraph*{Conjecture $\widehat{\ref{subloop.conjecture}}$:}
{\em
 Let $\Phi:{\mathcal{ A}}^{\mathbb{N}}{{\longrightarrow}}{\mathcal{ A}}^{\mathbb{N}}$ be a QGCA.
If $\mu$ is a $\Phi$-ergodic and $ {{{\boldsymbol{\sigma}}}^{}} $-invariant measure,
and $\entr{ {{{\boldsymbol{\sigma}}}^{}} }>0$, then $\mu$ is the uniform measure on ${\mathcal{ B}}^{\mathbb{N}}$,
for some ${\mathcal{ B}}\prec{\mathcal{ A}}$. }

 \medskip

  It is Conjecture $\widehat{\ref{subloop.conjecture}}$ which we'll refute in \S\ref{S:group}.

\section{Multiplication CA on Nonabelian Groups \label{S:group}}

Let ${\mathbb{N}}=\{0,1,2,3,\ldots\}$, and let $\mu$ be a
$ {{{\boldsymbol{\sigma}}}^{}} $-invariant measure on ${\mathcal{ A}}^{\mathbb{N}}$.  
Let ${\widetilde{\mathbb{N}}}=\{1,2,3,\ldots\}$
For any ${\mathbf{ a}} \in{\mathcal{ A}}^{{\widetilde{\mathbb{N}}}}$, let $\mu_{\mathbf{ a}}$ be the conditional measure induced
by ${\mathbf{ a}}$ on the zeroth coordinate.  That is, for any $b\in{\mathcal{ A}}$,
\[
  \mu_{\mathbf{ a}}(b) \quad=\quad \mu\left[x_0=b \ | \  {\mathbf{ x}}\raisebox{-0.3em}{$\left|_{{\widetilde{\mathbb{N}}}}\right.$} 
                                \ =  \ {\mathbf{ a}}\right].
\qquad\mbox{(where ${\mathbf{ x}}\in{\mathcal{ A}}^{\mathbb{N}}$ is a $\mu$-random sequence)}
\]
Let ${\widetilde{\mu}}$ be the projection of $\mu$ onto ${\mathcal{ A}}^{{\widetilde{\mathbb{N}}}}$.  Then we
have the following disintegration \cite{SchwartzDisintegrate}:
 \begin{equation} 
\label{disintegration}
  \mu \quad=\quad\int_{{\mathcal{ A}}^{{\widetilde{\mathbb{N}}}}} \left(\mu_{\mathbf{ a}} \otimes \delta_{\mathbf{ a}}\right) \ d{\widetilde{\mu}}[{\mathbf{ a}}].
 \end{equation} 
Suppose ${\mathcal{ A}}$ is a finite (possibly nonabelian) group, and let
${\mathcal{ C}}\prec{\mathcal{ A}}$ be a subgroup.  We call $\mu$ a {\bf ${\mathcal{ C}}$-measure} if, for
${{\forall}_{\mu}\;} {\mathbf{ a}}\in{\mathcal{ A}}^{\mathbb{N}}$, \ ${\sf supp}\left(\mu_{\mathbf{ a}}\right)$ is a right coset of ${\mathcal{ C}}$, and
$\mu_{\mathbf{ a}}$ is uniformly distributed on this coset.  It follows:

\begin{lemma}{\sf \label{B.measure.lemma}}  
 \setcounter{enumi}{\thethm} \begin{list}{{\bf (\alph{enumii})}}{\usecounter{enumii}} 			{\setlength{\leftmargin}{0em} 			\setlength{\rightmargin}{0em}} 
   \item If $\mu$ is a ${\mathcal{ C}}$-measure, then $h(\mu, {{{\boldsymbol{\sigma}}}^{}} ) \ = \ \log_2|{\mathcal{ C}}|$.

   \item Let $\{e\}$ be the identity subgroup.  Then
  $\left( \ \rule[-0.5em]{0em}{1em}       \begin{minipage}{40em}       \begin{tabbing}         $h(\mu, {{{\boldsymbol{\sigma}}}^{}} )  =  0$        \end{tabbing}      \end{minipage} \ \right)
\iff \left( \ \rule[-0.5em]{0em}{1em}       \begin{minipage}{40em}       \begin{tabbing}         $\mu$ is an $\{e\}$-measure        \end{tabbing}      \end{minipage} \ \right)$.

   \item  $\left( \ \rule[-0.5em]{0em}{1em}       \begin{minipage}{40em}       \begin{tabbing}         $\mu$ is an ${\mathcal{ A}}$-measure        \end{tabbing}      \end{minipage} \ \right) \iff
 \left( \ \rule[-0.5em]{0em}{1em}       \begin{minipage}{40em}       \begin{tabbing}         $\mu$ is the uniform measure on ${\mathcal{ A}}^{\mathbb{N}}$        \end{tabbing}      \end{minipage} \ \right)$.\hrulefill\ensuremath{\Box}
\end{list}
 \end{lemma}

 Let $\Phi:{\mathcal{ A}}^{\mathbb{N}}{{\longrightarrow}}{\mathcal{ A}}^{\mathbb{N}}$ be the {\bf nearest neighbour
multiplication} CA, having local map $\phi(a_0,a_1) = a_0\cdot a_1$.
This type of CA was previously studied in
\cite{MooreQuasi,PivatoNACA}. Our goal is to prove:

\begin{thm}{\sf \label{NN.mult.inv.meas}}  
  If $\mu$ is $ {{{\boldsymbol{\sigma}}}^{}} $-invariant and $\Phi$-ergodic,
then $\mu$ is a ${\mathcal{ C}}$-measure for some  ${\mathcal{ C}}\prec{\mathcal{ A}}$.\hrulefill\ensuremath{\Box}
 \end{thm}

  \medskip         \refstepcounter{thm} {\bf Example \thethm:}  \setcounter{enumi}{\thethm} \begin{list}{(\alph{enumii})}{\usecounter{enumii}} 			{\setlength{\leftmargin}{0em} 			\setlength{\rightmargin}{0em}}   
\item Let ${\mathcal{ C}}\prec{\mathcal{ A}}$ be any subgroup, and let $\mu$ be the uniform
measure on ${\mathcal{ C}}^{\mathbb{N}}$.  Then $\mu$ is a ${\mathcal{ C}}$-measure (for any ${\mathbf{ a}}\in{\mathcal{ A}}^{{\widetilde{\mathbb{N}}}}$, \ $\mu_{\mathbf{ a}}$ is uniform on ${\mathcal{ C}}$), and $\mu$ is
$ {{{\boldsymbol{\sigma}}}^{}} $-invariant and $\Phi$-ergodic.

\item \label{counterexample}
Let ${\mathcal{ Q}}=\{\pm1,\ \pm{\mathbf{ i}}, \ \pm{\mathbf{ j}}, \ \pm{\mathbf{ k}}\}$ be the Quaternion group
\cite[\S1.5]{DummitFoote},
and let $\Phi_{\scriptscriptstyle {\mathcal{ Q}}}:{\mathcal{ Q}}^{\mathbb{N}}{{\longrightarrow}}{\mathcal{ Q}}^{\mathbb{N}}$ be the nearest neighbour
multiplication CA.  It follows:
\[
\begin{array}{rrclcl}
 \mbox{If} & {\mathbf{ p}} &=& [{\mathbf{ i}},{\mathbf{ j}},{\mathbf{ k}},{\mathbf{ i}},{\mathbf{ j}},{\mathbf{ k}},{\mathbf{ i}},{\mathbf{ j}},{\mathbf{ k}},\ldots] \\ 
 \mbox{then} & \Phi_{\scriptscriptstyle {\mathcal{ Q}}}({\mathbf{ p}}) &=& [{\mathbf{ k}},{\mathbf{ i}},{\mathbf{ j}},{\mathbf{ k}},{\mathbf{ i}},{\mathbf{ j}},{\mathbf{ k}},{\mathbf{ i}},{\mathbf{ j}},\ldots], \\ 
 \mbox{and} & \Phi_{\scriptscriptstyle {\mathcal{ Q}}}^2({\mathbf{ p}}) &=& [{\mathbf{ j}},{\mathbf{ k}},{\mathbf{ i}},{\mathbf{ j}},{\mathbf{ k}},{\mathbf{ i}},{\mathbf{ j}},{\mathbf{ k}},{\mathbf{ i}},\ldots], \\ 
 \mbox{and} & \Phi_{\scriptscriptstyle {\mathcal{ Q}}}^3({\mathbf{ p}}) &=& [{\mathbf{ i}},{\mathbf{ j}},{\mathbf{ k}},{\mathbf{ i}},{\mathbf{ j}},{\mathbf{ k}},{\mathbf{ i}},{\mathbf{ j}},{\mathbf{ k}},\ldots]
&=&{\mathbf{ p}}.
\end{array}
\]
  Let $\mu_{\scriptscriptstyle {\mathcal{ Q}}}$ be the probability measure on ${\mathcal{ Q}}^{\mathbb{N}}$ assigning probability
$1/3$ to each of ${\mathbf{ p}}$, $\Phi_{\scriptscriptstyle {\mathcal{ Q}}}({\mathbf{ p}})$ and $\Phi_{\scriptscriptstyle {\mathcal{ Q}}}^2({\mathbf{ p}})$.  Then
$\mu_{\scriptscriptstyle {\mathcal{ Q}}}$ is $ {{{\boldsymbol{\sigma}}}^{}} $-invariant and $\Phi_{\scriptscriptstyle {\mathcal{ Q}}}$-ergodic.

Now, let ${\mathcal{ C}}$ be any other group, and let ${\mathcal{ A}} = {\mathcal{ C}} \times {\mathcal{ Q}}$.
Identify ${\mathcal{ C}}$ with ${\mathcal{ C}}\times\{1\}\ \prec \ {\mathcal{ A}}$; then ${\mathcal{ C}}$ is a
normal subgroup of ${\mathcal{ A}}$, and ${\mathcal{ Q}}={\mathcal{ A}}/{\mathcal{ C}}$.  The cosets of ${\mathcal{ C}}$ all
have the form ${\mathcal{ C}}\times\{q\}$ for some $q\in{\mathcal{ Q}}$.
  There is a natural identification 
${\mathcal{ A}}^{\mathbb{N}} \ \cong \  {\mathcal{ C}}^{\mathbb{N}} \times {\mathcal{ Q}}^{\mathbb{N}}$, given:
\[
  \left[\rule[-0.5em]{0em}{1em} (c_0,q_0), \ (c_1,q_1), \ (c_2,q_2), \ldots \right]
\quad \longleftrightarrow \quad
\left(\rule[-0.5em]{0em}{1em}  \left[c_0,c_1,c_2, \ldots \right]; \ \
 \left[ q_0,q_1q_2,\ldots \right] \right)
\]
 Let $\mu_{\scriptscriptstyle {\mathcal{ C}}}$ be the uniform Bernoulli measure on ${\mathcal{ C}}^{\mathbb{N}}$, and
let $\mu = \mu_{\scriptscriptstyle {\mathcal{ C}}} \otimes \mu_{\scriptscriptstyle {\mathcal{ Q}}}$.  

\setcounter{claimcount}{0}
\refstepcounter{claimcount}                {\bf Claim \theclaimcount: \ }{\sl  $\mu$ is a ${\mathcal{ C}}$-measure.}
\bprf
 Suppose ${\mathbf{ a}}\in{\mathcal{ A}}^{\mathbb{N}}$ is a $\mu$-random sequence.  Then 
$ {\mathbf{ a}} \ = \   \left( {\mathbf{ c}}, {\mathbf{ q}} \right)$,
where ${\mathbf{ q}} \in \left\{{\mathbf{ p}}, \ \Phi_{\scriptscriptstyle {\mathcal{ Q}}}({\mathbf{ p}}), \ \Phi_{\scriptscriptstyle {\mathcal{ Q}}}^2({\mathbf{ p}})\right\}$,
(with probability $1/3$ each), and ${\mathbf{ c}} = (c_0,c_1,c_2,\ldots)$ is a sequence of independent, uniformly distributed random elements of
${\mathcal{ C}}$.  The coordinates 
\[
 \left[a_1,a_2,a_3,\ldots \right]\quad =\quad 
\left[\rule[-0.5em]{0em}{1em} (c_1,q_1), \ (c_2,q_2), \ (c_2,q_2), \ldots \right]
\]
 determine ${\mathbf{ q}}$, and thus, determine $q_0$.  Thus, $\mu_{[a_1,a_2,a_3,\ldots]}$ is uniformly distributed on the coset ${\mathcal{ C}}\times\{q_0\}$.
 {\tt \dotfill~$\Diamond$~[Claim~\theclaimcount] }\end{list} 

  Let $\Phi:{\mathcal{ A}}^{\mathbb{N}}{{\longrightarrow}}{\mathcal{ A}}^{\mathbb{N}}$ be the nearest neighbour
multiplication map on ${\mathcal{ A}}^{\mathbb{N}}$.
 
\refstepcounter{claimcount}                {\bf Claim \theclaimcount: \ }{\sl  $\mu$ is $\Phi$-ergodic and $ {{{\boldsymbol{\sigma}}}^{}} $-invariant.}
\bprf   
$\mu$ is clearly $ {{{\boldsymbol{\sigma}}}^{}} $-invariant.

  {\em $\mu$ is $\Phi$-invariant}: \quad
  Let $\Phi_{\scriptscriptstyle {\mathcal{ C}}}:{\mathcal{ C}}^{\mathbb{N}}{{\longrightarrow}}{\mathcal{ C}}^{\mathbb{N}}$ be the nearest neighbour
multiplication map on ${\mathcal{ C}}^{\mathbb{N}}$.  Then $\Phi \ = \ \Phi_{\scriptscriptstyle {\mathcal{ C}}} \times \Phi_{\scriptscriptstyle {\mathcal{ Q}}}$.  Thus,
$\Phi(\mu) \ = \ \Phi_{\scriptscriptstyle {\mathcal{ C}}}(\mu_{\scriptscriptstyle {\mathcal{ C}}}) \otimes \Phi_{\scriptscriptstyle {\mathcal{ Q}}}(\mu_{\scriptscriptstyle {\mathcal{ Q}}})
\ = \ \mu_{\scriptscriptstyle {\mathcal{ C}}} \otimes \mu_{\scriptscriptstyle {\mathcal{ Q}}} \ = \ \mu$. 

 {\em $\mu$ is $\Phi$-ergodic}: \quad The system $({\mathcal{ C}},\Phi_{\scriptscriptstyle {\mathcal{ C}}},\mu_{\scriptscriptstyle {\mathcal{ C}}})$ is
mixing \cite[Thm 6.3]{Kleveland}, thus weakly mixing.
The system $({\mathcal{ Q}},\Phi_{\scriptscriptstyle {\mathcal{ Q}}},\mu_{\scriptscriptstyle {\mathcal{ Q}}})$ is
ergodic.  Thus, the product system $({\mathcal{ A}},\Phi,\mu) \ = \ \left({\mathcal{ C}}\times{\mathcal{ Q}},
\ \Phi_{\scriptscriptstyle {\mathcal{ C}}}\times\Phi_{\scriptscriptstyle {\mathcal{ Q}}}, \ \mu_{\scriptscriptstyle {\mathcal{ C}}}\otimes\mu_{\scriptscriptstyle {\mathcal{ Q}}}\right)$ is also ergodic \cite[Thm. 2.6.1]{Petersen}.
 {\tt \dotfill~$\Diamond$~[Claim~\theclaimcount] }\end{list}
  Observe that $h(\mu, {{{\boldsymbol{\sigma}}}^{}} ) \ = \ h(\mu_{\scriptscriptstyle {\mathcal{ C}}}, {{{\boldsymbol{\sigma}}}^{}} ) \ = \ \log_2|{\mathcal{ C}}|$.
Thus, $\mu$ is a $\Phi$-ergodic, $ {{{\boldsymbol{\sigma}}}^{}} $-invariant measure
of nonzero entropy.  However, ${\sf supp}\left(\mu\right) \neq {\mathcal{ B}}^{\mathbb{N}}$ for any subgroup
${\mathcal{ B}}\prec{\mathcal{ A}}$.  This contradicts Conjecture 
$\widehat{\ref{subloop.conjecture}}$. 	\hrulefill\end{list}   			

\begin{cor}{\sf }  
  Let $h_{\rm max} \ = \ \max{\left\{ \log_2|{\mathcal{ C}}| \; ; \; {\mathcal{ C}} \ \mbox{a proper subgroup of} \ {\mathcal{ A}} \right\} }$.
(In particular, if 
${\mathcal{ A}}$ has no nontrivial proper subgroups, then $h_{\rm max}=0$.)

    If $\mu$ is  $ {{{\boldsymbol{\sigma}}}^{}} $-invariant and $\Phi$-ergodic,
and $h(\mu, {{{\boldsymbol{\sigma}}}^{}} ) \ >  \ h_{\rm max}$, then $\mu$ is the uniform measure. \end{cor}
\bprf
  Theorem \ref{NN.mult.inv.meas} says $\mu$ must be a ${\mathcal{ C}}$-measure
for some subgroup ${\mathcal{ C}}\prec{\mathcal{ A}}$.  But if ${\mathcal{ B}}$ is any proper subgroup, then
 $h(\mu, {{{\boldsymbol{\sigma}}}^{}} )  \ > \  h_{\rm max} \ \geq  \ \log_2|{\mathcal{ B}}|$,
so Lemma \ref{B.measure.lemma}{\bf(a)} says ${\mathcal{ C}}$ can't be ${\mathcal{ B}}$.
Thus, ${\mathcal{ C}}={\mathcal{ A}}$.  Then Lemma \ref{B.measure.lemma}{\bf(c)} says that $\mu$ is
the uniform measure.
 {\tt \hrulefill $\Box$ } \end{list}  \medskip  

        \refstepcounter{thm}                     \begin{list}{} 			{\setlength{\leftmargin}{1em} 			\setlength{\rightmargin}{0em}}        \item {\bf Example \thethm:} If $p$ and $q$ are prime and $p$ divides $q-1$, then there
is a unique nonabelian group of order $pq$ \cite[\S5.5]{DummitFoote}.
For example, let $p=3$ and $q=7$ and let ${\mathcal{ A}}$ be the unique
nonabelian group of order $21$.  Then $h_{\rm max}=\log_2(7)\approx 2.807 < 
4.392 \approx \log_2(21)$.  Hence, if
 $\mu$ is $ {{{\boldsymbol{\sigma}}}^{}} $-invariant and $\Phi$-ergodic,
and $h(\mu, {{{\boldsymbol{\sigma}}}^{}} ) \geq 2.81$, then $\mu$ is the uniform measure.                   \hrulefill  \end{list}   			

  If $b\in{\mathcal{ A}}$, then we define {\bf (left) scalar multiplication} by
$b$ upon ${\mathcal{ A}}^{\mathbb{N}}$ in the obvious way: if
 ${\mathbf{ c}}=[c_0,c_1,c_2,\ldots] \in
{\mathcal{ A}}^{\mathbb{N}}$, then $b\cdot{\mathbf{ c}} \ = \ [bc_0, \ bc_1, \ bc_2,\ldots]$.
  For any sequence ${\mathbf{ a}}=[a_1,a_2,a_3,\ldots]$ in ${\mathcal{ A}}^{{\widetilde{\mathbb{N}}}}$ and any
$b\in{\mathcal{ A}}$, let $[b,{\mathbf{ a}}]$ denote the sequence $[b,a_1,a_2,a_3,\ldots]$ in
${\mathcal{ A}}^{\mathbb{N}}$.  Recall the conjugacy $\Xi:{\mathcal{ A}}^{\mathbb{N}}{{\longrightarrow}}{\mathcal{ A}}^{\mathbb{N}}$
and the dual cellular automaton ${\widehat{\Phi }}:{\mathcal{ A}}^{\mathbb{N}}{{\longrightarrow}}{\mathcal{ A}}^{\mathbb{N}}$
introduced in \S\ref{S:dualCA}.

\begin{lemma}{\sf \label{Xi.scalar}}  
Let ${\mathbf{ a}}\in{\mathcal{ A}}^{{\widetilde{\mathbb{N}}}}$, and suppose $\Xi({\mathbf{ a}}) \ = \ [b_0, b_1, b_2, \ldots]$.
  Then:
 \setcounter{enumi}{\thethm} \begin{list}{{\bf (\alph{enumii})}}{\usecounter{enumii}} 			{\setlength{\leftmargin}{0em} 			\setlength{\rightmargin}{0em}}
  \item $\Xi[e,{\mathbf{ a}}] \ = \ \left[e, \ b_0, \ b_0b_1, \ b_0 b_1 b_2, \ b_0 b_1 b_2 b_3, \ldots\right]$,
  \item For any $b\in{\mathcal{ A}}$, \quad $\Xi[b,{\mathbf{ a}}] \ = \ b\cdot\Xi[e,{\mathbf{ a}}]$.\hrulefill\ensuremath{\Box}
\end{list}
 \end{lemma}

  Say that an element ${\mathbf{ g}}\in{\mathcal{ A}}^{\mathbb{N}}$ is {\bf $(\Phi,\mu)$-generic}
if, for any cylinder set ${\mathbf{ U}}\subset{\mathcal{ A}}^{\mathbb{N}}$,
\[
  \mu[{\mathbf{ U}}] \quad=\quad \lim_{N{\rightarrow}{\infty}} \frac{1}{N} \sum_{n=1}^N {{{\mathsf{ 1\!\!1}}}_{{{\mathbf{ U}}}}}\left(\Phi^n({\mathbf{ g}})\right).
\]
  The Birkhoff Ergodic Theorem says that $\mu$-almost all points in
${\mathcal{ A}}^{\mathbb{N}}$ are $(\Phi,\mu)$-generic.

  Let $\nu = \Xi(\mu)$.  It follows that $\nu$ is ${\widehat{\Phi }}$-invariant
and $ {{{\boldsymbol{\sigma}}}^{}} $-ergodic.   Lemma \ref{phi.to.shift.conj} implies:

\begin{lemma}{\sf \label{gen.gen}}  
 Let ${\mathbf{ g}}\in{\mathcal{ A}}^{\mathbb{N}}$.  Then 
$\left( \ \rule[-0.5em]{0em}{1em}       \begin{minipage}{40em}       \begin{tabbing}         ${\mathbf{ g}}$ is $(\Phi,\mu)$-generic        \end{tabbing}      \end{minipage} \ \right)\iff
\left( \ \rule[-0.5em]{0em}{1em}       \begin{minipage}{40em}       \begin{tabbing}          $\Xi({\mathbf{ g}})$ is $( {{{\boldsymbol{\sigma}}}^{}} ,\nu)$-generic        \end{tabbing}      \end{minipage} \ \right).$\hrulefill\ensuremath{\Box}
 \end{lemma}

\begin{lemma}{\sf \label{foobaz1}}  
 Let ${\mathbf{ a}}\in{\mathcal{ A}}^{{\widetilde{\mathbb{N}}}}$, and let $b,b'\in{\mathcal{ A}}$.  Suppose both
$[b,{\mathbf{ a}}]$ and $[b',{\mathbf{ a}}]$ are $(\Phi,\mu)$-generic. 

 If $c = b'\cdot b^{-1} $,  then $\nu$ is invariant under (left) scalar
 multiplication by $c$. 
In other words, for any measurable subset ${\mathbf{ U}}\subset{\mathcal{ A}}^{\mathbb{N}}$,
\quad $\mu[c\cdot{\mathbf{ U}}] \ = \ \mu[{\mathbf{ U}}]$.
 \end{lemma}
\bprf
  Let ${\mathbf{ g}} = \Xi[b,{\mathbf{ a}}]$ and ${\mathbf{ g}}' = \Xi[b',{\mathbf{ a}}]$.  Then
 $b' = c \cdot b$, so
 \begin{equation} 
\label{foobar}
  {\mathbf{ g}}' \quad=\quad\Xi[b',{\mathbf{ a}}] 
\quad\displaystyle\raisebox{-1ex}{$\overline{\overline{{\scriptscriptstyle{\mathrm{(L\ref{Xi.scalar})}}}}}$}\quad b'\cdot \Xi[e,{\mathbf{ a}}] 
\quad=\quad cb\cdot \Xi[e,{\mathbf{ a}}] 
\quad\displaystyle\raisebox{-1ex}{$\overline{\overline{{\scriptscriptstyle{\mathrm{(L\ref{Xi.scalar})}}}}}$}\quad c\cdot \Xi[b,{\mathbf{ a}}] 
\quad=\quad c\cdot {\mathbf{ g}},
 \end{equation} 
  where {\bf(L\ref{Xi.scalar})} is by Lemma \ref{Xi.scalar}{\bf(b)}.
Next, Lemma \ref{gen.gen} says that ${\mathbf{ g}}$ and ${\mathbf{ g}}'$ are both 
 $( {{{\boldsymbol{\sigma}}}^{}} ,\nu)$-generic.  Thus, for any cylinder set ${\mathbf{ U}}\subset{\mathcal{ A}}^{\mathbb{N}}$, 
\begin{eqnarray*}
  \nu[{\mathbf{ U}}] &\displaystyle\raisebox{-1ex}{$\overline{\overline{{\scriptscriptstyle{\mathrm{(g1)}}}}}$}&
 \lim_{N{\rightarrow}{\infty}} \frac{1}{N} \sum_{n=1}^N {{{\mathsf{ 1\!\!1}}}_{{{\mathbf{ U}}}}}\left( {{{\boldsymbol{\sigma}}}^{n}} ({\mathbf{ g}})\right)
\quad\displaystyle\raisebox{-1ex}{$\overline{\overline{{\scriptscriptstyle{\mathrm{(eq\ref{foobar})}}}}}$}\quad 
\lim_{N{\rightarrow}{\infty}} \frac{1}{N} \sum_{n=1}^N {{{\mathsf{ 1\!\!1}}}_{{{\mathbf{ U}}}}}\left( {{{\boldsymbol{\sigma}}}^{n}} (c^{-1}\cdot {\mathbf{ g}}')\right)
\\&=&
\lim_{N{\rightarrow}{\infty}} \frac{1}{N} \sum_{n=1}^N {{{\mathsf{ 1\!\!1}}}_{{(c\cdot{\mathbf{ U}})}}}\left( {{{\boldsymbol{\sigma}}}^{n}} ({\mathbf{ g}}')\right)
\quad\displaystyle\raisebox{-1ex}{$\overline{\overline{{\scriptscriptstyle{\mathrm{(g2)}}}}}$}\quad
\nu[c\cdot {\mathbf{ U}}].
\end{eqnarray*}
 {\bf(g1)} is because ${\mathbf{ g}}$ is generic,
{\bf(eq\ref{foobar})} is by eqn. (\ref{foobar}), and
{\bf(g2)} is because ${\mathbf{ g}}'$ is generic.
 {\tt \hrulefill $\Box$ } \end{list}  \medskip   

  We next show that the hypothesis of Lemma \ref{foobaz1} is
not vacuous.  Let
\[
 {\widetilde{\mathbf{ F}}}_2\quad=\quad{\left\{ {\mathbf{ a}}\in{\mathcal{ A}}^{{\widetilde{\mathbb{N}}}} \; ; \; {\sf card}\left[{\sf supp}\left(\mu_{\mathbf{ a}}\right)\right]\geq 2 \right\} }.
\]
\begin{lemma}{\sf \label{foobaz0}}   If $h(\mu, {{{\boldsymbol{\sigma}}}^{}} )>0$, then ${\widetilde{\mu}}[{\widetilde{\mathbf{ F}}}_2]  \ > \ 0$. \end{lemma}
\bprf
Let ${\widetilde{\mathbf{ F}}}_1={\left\{ {\mathbf{ a}}\in{\mathcal{ A}}^{{\widetilde{\mathbb{N}}}} \; ; \; {\sf card}\left[{\sf supp}\left(\mu_{\mathbf{ a}}\right)\right]\geq 1 \right\} }$.

\refstepcounter{claimcount}                {\bf Claim \theclaimcount: \ }{\sl  ${\widetilde{\mu}}[{\widetilde{\mathbf{ F}}}_1]=1$.}
\bprf
$\displaystyle 1 \ = \quad  \mu[{\mathcal{ A}}^{\mathbb{N}}] 
\quad \displaystyle\raisebox{-1ex}{$\overline{\overline{{\scriptscriptstyle{\mathrm{(eq\ref{disintegration})}}}}}$} \quad  \int_{{\mathcal{ A}}^{{\widetilde{\mathbb{N}}}}} \mu_{\mathbf{ a}}[{\mathcal{ A}}] \ d{\widetilde{\mu}}[{\mathbf{ a}}]
\ = \  \int_{{\widetilde{\mathbf{ F}}}_1} 1\ d{\widetilde{\mu}}[{\mathbf{ a}}]
\ = \ {\widetilde{\mu}}[{\widetilde{\mathbf{ F}}}_1]$. 
 {\tt \dotfill~$\Diamond$~[Claim~\theclaimcount] }\end{list}
If $\rho$ is a measure on ${\mathcal{ A}}$, define $\displaystyle H(\rho) =  - \sum_{b\in{\mathcal{ A}}}
\rho\{b\}\log_2\left(\rho\{b\}\right)$. Recall \cite[Prop. 5.2.12]{Petersen} that 
 \begin{equation} 
\label{entropy.int}
h(\mu, {{{\boldsymbol{\sigma}}}^{}} ) \ = \ \int_{{\mathcal{ A}}^{{\widetilde{\mathbb{N}}}}} H(\mu_{\mathbf{ a}}) \ d{\widetilde{\mu}}[{\mathbf{ a}}].
 \end{equation}  
 
\refstepcounter{claimcount}                {\bf Claim \theclaimcount: \ }{\sl  If ${\widetilde{\mu}}[{\widetilde{\mathbf{ F}}}_2] = 0$, then $H(\mu_{\mathbf{ a}}) \ = \ 0$ for $\forall_{{\widetilde{\mu}}}
\ {\mathbf{ a}} \in {\mathcal{ A}}^{{\widetilde{\mathbb{N}}}}$.}
\bprf
 Let ${\widetilde{\mathbf{ F}}}_* \ = \ {\widetilde{\mathbf{ F}}}_1 \setminus {\widetilde{\mathbf{ F}}}_2
\ = \ 
{\left\{ {\mathbf{ a}}\in{\mathcal{ A}}^{{\widetilde{\mathbb{N}}}} \; ; \; {\sf card}\left[{\sf supp}\left(\mu_{\mathbf{ a}}\right)\right]= 1 \right\} }$. \
If ${\widetilde{\mu}}[{\widetilde{\mathbf{ F}}}_2] = 0$, then $\mu[{\widetilde{\mathbf{ F}}}_*] = \mu[{\widetilde{\mathbf{ F}}}_1] -
\mu[{\widetilde{\mathbf{ F}}}_2] = 1$.  Thus, there is measurable function
$\gamma:{\mathcal{ A}}^{{\widetilde{\mathbb{N}}}}{{\longrightarrow}}{\mathcal{ A}}$ so that
$\mu_{\mathbf{ a}}(\gamma({\mathbf{ a}})) \ = \ 1$ for $\forall_{{\widetilde{\mu}}} \ {\mathbf{ a}}\in{\mathcal{ A}}^{{\widetilde{\mathbb{N}}}}$.
 Hence, $H(\mu_{\mathbf{ a}})=0$, for $\forall_{{\widetilde{\mu}}} \ {\mathbf{ a}}\in{\mathcal{ A}}^{{\widetilde{\mathbb{N}}}}$.
 {\tt \dotfill~$\Diamond$~[Claim~\theclaimcount] }\end{list}
  Claim 2 and  equation (\ref{entropy.int}) imply
that $h(\mu, {{{\boldsymbol{\sigma}}}^{}} )=0$, contradicting our hypothesis. 
 {\tt \hrulefill $\Box$ } \end{list}  \medskip  

  Let ${\widetilde{\mathbf{ G}}} \ = \ {\left\{ {\mathbf{ a}}\in{\mathcal{ A}}^{{\widetilde{\mathbb{N}}}} \; ; \; \mbox{$[b,{\mathbf{ a}}]$ is $(\Phi,\mu)$-generic for every $b\in{\sf supp}\left(\mu_{\mathbf{ a}}\right)$} \right\} }$.

\begin{lemma}{\sf \label{foobaz3}}    ${\widetilde{\mu}}[{\widetilde{\mathbf{ G}}}]=1$. \end{lemma}
\bprf
  Suppose not.   Let
\[{\widetilde{\mathbf{ H}}} 
\ = \ {\mathcal{ A}}^{{\widetilde{\mathbb{N}}}}\setminus{\widetilde{\mathbf{ G}}} \ = \ 
{\left\{ {\mathbf{ a}}\in{\mathcal{ A}}^{{\widetilde{\mathbb{N}}}} \; ; \; \mbox{$[b,{\mathbf{ a}}]$ is {\em not} 
$(\Phi,\mu)$-generic for some $b\in{\sf supp}\left(\mu_{\mathbf{ a}}\right)$} \right\} }.
\]
  For every ${\mathbf{ h}}\in{\widetilde{\mathbf{ H}}}$, let ${\mathcal{ B}}_{\mathbf{ h}} \ = \ {\left\{  b\in{\sf supp}\left(\mu_{\mathbf{ h}}\right) \; ; \; \mbox{$[b,{\mathbf{ a}}]$ is {\em not} $(\Phi,\mu)$-generic.} \right\} }$.  Thus,
 \begin{equation} 
\label{foobaz3.1}
\forall\ {\mathbf{ h}}\in{\widetilde{\mathbf{ H}}}, \qquad \mu_{\mathbf{ h}}[{\mathcal{ B}}_{\mathbf{ h}}]\quad>\quad 0.
 \end{equation} 
Let ${\mathbf{ H}} \ = \ {\left\{ [b,{\mathbf{ h}}] \; ; \; {\mathbf{ h}}\in{\widetilde{\mathbf{ H}}}, \ b\in{\mathcal{ B}}_{\mathbf{ h}}  \right\} }$. 
If ${\widetilde{\mu}}[{\widetilde{\mathbf{ G}}}]<1$, then ${\widetilde{\mu}}[{\widetilde{\mathbf{ H}}}]>0$.  Thus,
\[
  \mu[{\mathbf{ H}}] \quad\displaystyle\raisebox{-1ex}{$\overline{\overline{{\scriptscriptstyle{\mathrm{(eq\ref{disintegration})}}}}}$}\quad \int_{{\widetilde{\mathbf{ H}}}} \mu_{\mathbf{ h}}[{\mathcal{ B}}_{\mathbf{ h}}] \ d{\widetilde{\mu}}[{\mathbf{ h}}]
\quad \raisebox{-1ex}{${{\displaystyle >} \atop {\scriptscriptstyle{\mathrm{(eq\ref{foobaz3.1})}}}}$} \quad 0,
\]
  Now let ${\mathbf{ G}} \ = \ 
{\left\{ {\mathbf{ g}}\in{\mathcal{ A}}^{{\mathbb{N}}} \; ; \; \mbox{${\mathbf{ g}}$ is $(\Phi,\mu)$-generic} \right\} }$.
Then the Birkhoff Ergodic Theorem says $\mu[{\mathbf{ G}}]=1$.  But
clearly ${\mathbf{ G}}\subset{\mathcal{ A}}^{\mathbb{N}}\setminus{\mathbf{ H}}$, so if $\mu[{\mathbf{ H}}]>0$,
then $\mu[{\mathbf{ G}}]<1$.  Contradiction.
 {\tt \hrulefill $\Box$ } \end{list}  \medskip  

Let ${\widetilde{\mathbf{ I}}}={\left\{ {\mathbf{ a}}\in{\mathcal{ A}}^{{\widetilde{\mathbb{N}}}} \; ; \; \mbox{there are distinct $b,b'\in{\mathcal{ A}}$
so that $[b,{\mathbf{ a}}]$ and $[b',{\mathbf{ a}}]$ are $(\Phi,\mu)$-generic} \right\} }$.

\begin{lemma}{\sf \label{foobar30}}  
If $h(\mu, {{{\boldsymbol{\sigma}}}^{}} )>0$, then  ${\widetilde{\mu}}[{\widetilde{\mathbf{ I}}}]>0$
(so the hypothesis of Lemma \ref{foobaz1} is nonvacuous).
 \end{lemma}
\bprf
  Observe that ${\widetilde{\mathbf{ I}}}\ \supset \ {\widetilde{\mathbf{ F}}}_2\cap{\widetilde{\mathbf{ G}}}$.  Now combine Lemmas \ref{foobaz0} and \ref{foobaz3}.
 {\tt \hrulefill $\Box$ } \end{list}  \medskip

\begin{lemma}{\sf \label{foobaz2}}  
If $\nu$ is invariant under
scalar multiplication by $c$,  then, for $\forall_{{\widetilde{\mu}}} \ 
{\mathbf{ a}}\in{\mathcal{ A}}^{\widetilde{\mathbb{N}}}$,\quad $\mu_{\mathbf{ a}}$ is invariant under 
left multiplication by $c$.
 \end{lemma}
\bprf  Let $b\in{\mathcal{ A}}$ and let $b' = c\cdot b$.  Define
$\beta:{\mathcal{ A}}^{{\widetilde{\mathbb{N}}}}{{\longrightarrow}}{\mathbb{R}}$ by $\beta({\mathbf{ a}}) \ = \ \mu_{\mathbf{ a}}(b)$.
Likewise, let $\beta'({\mathbf{ a}}) \ = \ \mu_{\mathbf{ a}}(b')$.  Then $\beta$ and
$\beta'$ are measurable, and we want to show that $\beta = \beta'$, \ 
${\widetilde{\mu}}$-\ae.  

  Define $\gamma:{\mathcal{ A}}^{\mathbb{N}}{{\longrightarrow}}{\mathcal{ A}}^{\mathbb{N}}$ by $\gamma[a_0,a_1,a_2,\ldots]
\ = \ [c\cdot a_0, \ a_1, \ a_2, \ldots]$.  Define $\Gamma:{\mathcal{ A}}^{\mathbb{N}}{{\longrightarrow}}{\mathcal{ A}}^{\mathbb{N}}$ by $\Gamma[{\mathbf{ a}}] = c\cdot {\mathbf{ a}}$.  

\refstepcounter{claimcount}                {\bf Claim \theclaimcount: \ }{\sl  \label{foobaz2.1}$\Xi\circ \gamma \ = \ \Gamma \circ \Xi$.}
\bprf 
 Generalize the reasoning behind
equation (\ref{foobar}).
 {\tt \dotfill~$\Diamond$~[Claim~\theclaimcount] }\end{list}

\refstepcounter{claimcount}                {\bf Claim \theclaimcount: \ }{\sl  \label{foobaz2.2}$\mu$ is $\gamma$-invariant.}
\bprf
 For any measurable subset ${\mathbf{ U}}\subset{\mathcal{ A}}^{\mathbb{N}}$,
\[
 \mu\left[\rule[-0.5em]{0em}{1em} \gamma({\mathbf{ U}}) \right] \quad \displaystyle\raisebox{-1ex}{$\overline{\overline{{\scriptscriptstyle{\mathrm{(D)}}}}}$}\quad
 \nu\left[\rule[-0.5em]{0em}{1em}\Xi\circ\gamma({\mathbf{ U}}) \right] \quad \displaystyle\raisebox{-1ex}{$\overline{\overline{{\scriptscriptstyle{\mathrm{(C\ref{foobaz2.1})}}}}}$}\quad
 \nu\left[\rule[-0.5em]{0em}{1em}\Gamma\circ\Xi({\mathbf{ U}}) \right] \quad \displaystyle\raisebox{-1ex}{$\overline{\overline{{\scriptscriptstyle{\mathrm{(I)}}}}}$}\quad
 \nu\left[\rule[-0.5em]{0em}{1em}\Xi({\mathbf{ U}}) \right] \quad \displaystyle\raisebox{-1ex}{$\overline{\overline{{\scriptscriptstyle{\mathrm{(D)}}}}}$}\quad
\mu\left[{\mathbf{ U}}\right].
\]
 {\bf(D)} is by definition of $\nu$.   {\bf(C\ref{foobaz2.1})} is Claim
\ref{foobaz2.1}.   {\bf(I)} is because  $\nu$ is $\Gamma$-invariant.
 {\tt \dotfill~$\Diamond$~[Claim~\theclaimcount] }\end{list}
  
\refstepcounter{claimcount}                {\bf Claim \theclaimcount: \ }{\sl  \label{foobaz2.3}
For any measurable subset ${\widetilde{\mathbf{ W}}} \subset {\mathcal{ A}}^{{\widetilde{\mathbb{N}}}}$,
\quad $\displaystyle
\int_{{\widetilde{\mathbf{ W}}}} \beta({\mathbf{ w}}) \ d{\widetilde{\mu}}[{\mathbf{ w}}]
\quad=\quad
\int_{{\widetilde{\mathbf{ W}}}} \beta'({\mathbf{ w}}) \ d{\widetilde{\mu}}[{\mathbf{ w}}]$.
}
\bprf
Let ${\mathbf{ U}} \ = \ [b]\times {\widetilde{\mathbf{ W}}}$, and let  ${\mathbf{ U}}' \ = \ \gamma({\mathbf{ U}})
 \ = \  [cb]\times {\widetilde{\mathbf{ W}}} \ = \ [b']\times {\widetilde{\mathbf{ W}}}$.  Then:
\[
\int_{{\widetilde{\mathbf{ W}}}} \beta({\mathbf{ w}}) \ d{\widetilde{\mu}}[{\mathbf{ w}}]
\quad\displaystyle\raisebox{-1ex}{$\overline{\overline{{\scriptscriptstyle{\mathrm{(eq\ref{disintegration})}}}}}$}\quad
\mu[{\mathbf{ U}}] \quad\displaystyle\raisebox{-1ex}{$\overline{\overline{{\scriptscriptstyle{\mathrm{(c\ref{foobaz2.2})}}}}}$}\quad 
\mu[{\mathbf{ U}}'] \quad\displaystyle\raisebox{-1ex}{$\overline{\overline{{\scriptscriptstyle{\mathrm{(eq\ref{disintegration})}}}}}$}\quad
 \int_{{\widetilde{\mathbf{ W}}}} \beta'({\mathbf{ w}}) \ d{\widetilde{\mu}}[{\mathbf{ w}}],
\]
 where {\bf(eq\ref{disintegration})} is by equation (\ref{disintegration}), and
  {\bf(c\ref{foobaz2.2})} is by Claim \ref{foobaz2.2}.
 {\tt \dotfill~$\Diamond$~[Claim~\theclaimcount] }\end{list}
  It follows from Claim \ref{foobaz2.3} that
 $\beta \ = \ \beta'$, ${\widetilde{\mu}}$-\ae.  
 {\tt \hrulefill $\Box$ } \end{list}  \medskip  

\bprf[Proof of Theorem \ref{NN.mult.inv.meas}]
\setcounter{claimcount}{0}
 Let ${\mathcal{ C}}$ be the set of all $c\in{\mathcal{ A}}$ so that there is some
${\mathbf{ a}}\in{\mathcal{ A}}^{{\widetilde{\mathbb{N}}}}$ and $b\in{\mathcal{ A}}$ with both $[b,{\mathbf{ a}}]$ and
$[(cb),{\mathbf{ a}}]$ being $(\Phi,\mu)$-generic. 

If $h(\mu, {{{\boldsymbol{\sigma}}}^{}} )=0$, then $\mu$ is an $\{e\}$-measure by Lemma
\ref{B.measure.lemma}(b).  So, assume $h(\mu, {{{\boldsymbol{\sigma}}}^{}} )\neq 0$; then
Lemma \ref{foobar30} implies that ${\mathcal{ C}}$ is nontrivial.

\refstepcounter{claimcount}                {\bf Claim \theclaimcount: \ }{\sl  \label{NN.mult.inv.meas.C1} ${\mathcal{ C}}$ is a group, and  $\mu_{\mathbf{ a}}$ is
invariant under (left) ${\mathcal{ C}}$-multiplication for
$\forall_{{\widetilde{\mu}}} \ {\mathbf{ a}}\in{\mathcal{ A}}^{{\widetilde{\mathbb{N}}}}$. }
\bprf 
 Lemma \ref{foobaz1}
says that $\nu$ is invariant under ${\mathcal{ C}}$-scalar multiplication.
 Let ${\mathcal{ D}}$ be the group generated by ${\mathcal{ C}}$.  Then ${\mathcal{ C}}\subseteq{\mathcal{ D}}$,
and $\nu$ is also invariant under ${\mathcal{ D}}$-scalar
multiplication.  Lemma \ref{foobaz2} implies that $\mu_{\mathbf{ a}}$ is
invariant under (left) ${\mathcal{ D}}$-multiplication for
$\forall_{{\widetilde{\mu}}} \ {\mathbf{ a}}\in{\mathcal{ A}}^{{\widetilde{\mathbb{N}}}}$.  It follows from Lemma 
\ref{foobaz3} that ${\mathcal{ D}}\subseteq{\mathcal{ C}}$, and hence, ${\mathcal{ C}}={\mathcal{ D}}$.
 {\tt \dotfill~$\Diamond$~[Claim~\theclaimcount] }\end{list}

\refstepcounter{claimcount}                {\bf Claim \theclaimcount: \ }{\sl  \label{NN.mult.inv.meas.C2}
For $\forall_{{\widetilde{\mu}}} \ {\mathbf{ a}}\in{\mathcal{ A}}^{{\widetilde{\mathbb{N}}}}$,\quad
${\sf supp}\left(\mu_{\mathbf{ a}}\right)$ is a (right) coset of ${\mathcal{ C}}$.}
\bprf
 For $\forall_{{\widetilde{\mu}}} \ {\mathbf{ a}}\in{\mathcal{ A}}^{{\widetilde{\mathbb{N}}}}$, Claim
\ref{NN.mult.inv.meas.C1} implies that ${\sf supp}\left(\mu_{\mathbf{ a}}\right)$ is a disjoint
union of cosets of ${\mathcal{ C}}$, and that $\mu_{\mathbf{ a}}$ is uniformly distributed
on each of these cosets.  Let
\[
  {\widetilde{\mathbf{ M}}} \quad=\quad{\left\{  {\mathbf{ a}}\in{\mathcal{ A}}^{{\widetilde{\mathbb{N}}}} \; ; \; {\sf supp}\left(\mu_{\mathbf{ a}}\right) \ 
\mbox{contains more than one coset of ${\mathcal{ C}}$} \right\} }.
\]
  We claim that ${\widetilde{\mu}}[{\widetilde{\mathbf{ M}}}] \ = \ 0$.  Suppose not.  Then Lemma
\ref{foobaz3} implies that $\mu[{\widetilde{\mathbf{ M}}}\cap{\widetilde{\mathbf{ G}}}]>0$.  So let
${\mathbf{ m}}\in{\widetilde{\mathbf{ M}}}\cap{\widetilde{\mathbf{ G}}}$, and find elements $b,b'\in{\sf supp}\left(\mu_{\mathbf{ m}}\right)$
living in {\em different} cosets, such that $[b,{\mathbf{ m}}]$ and $[b',{\mathbf{ m}}]$
are both $(\Phi,\mu)$-generic.  If $c = b^{-1}b'$, then $b' = cb$, so
$c\in{\mathcal{ C}}$.  But $b$ and $b'$ are in different cosets of ${\mathcal{ C}}$; hence,
$c\not\in{\mathcal{ C}}$.  Contradiction.
\hrulefill $\Diamond$ {\tt[Claim~\theclaimcount]~~$\Box$} \end{list} \end{list}  \medskip  
\section{\label{S:degree} Degree of QGCA relative to invariant measures}

  If $\mu$ is a $\Phi$-invariant measure, then $\Phi$ is {\bf $K$-to-1
  {\rm ($\mu$-\ae)}} if there is a subset ${\mathcal{ U}}\subset{{\mathcal{ A}}^{\mathbb{Z}}}$ such that:

\medskip

 1. $\mu[{\mathcal{ U}}]=1$.
\medskip

 2. $\Phi^{-1}({\mathcal{ U}}) \ \displaystyle\raisebox{-0.6ex}{$\overline{\overline{{\scriptstyle{\mathrm{\,\mu\,}}}}}$} \ {\mathcal{ U}}$.

 3. $\mu$-almost every element ${\mathbf{ u}}\in{\mathcal{ U}}$ has exactly $K$ preimages in ${\mathcal{ U}}$ ---ie.  $\left|\rule[-0.5em]{0em}{1em}{\mathcal{ U}}\cap\Phi^{-1}\{{\mathbf{ u}}\}\right| \ = \ K$.

\medskip

  We will generalize the methods of \cite{HostMaassMartinez} to prove:

\begin{thm}{\sf \label{loop.ca.K.2.1}}  
 Let $\Phi:{{\mathcal{ A}}^{\mathbb{Z}}}{{\longrightarrow}}{{\mathcal{ A}}^{\mathbb{Z}}}$ be a quasigroup CA, and let
$\mu$ be a measure which is $\Phi$-invariant and $ {{{\boldsymbol{\sigma}}}^{}} $-ergodic.
Let $|{\mathcal{ A}}| = N$.  Then there is some $K\in{\left[ 1..N \right]}$ so that
 \setcounter{enumi}{\thethm} \begin{list}{{\bf (\alph{enumii})}}{\usecounter{enumii}} 			{\setlength{\leftmargin}{0em} 			\setlength{\rightmargin}{0em}}
  \item $\entr{\Phi} \ = \ \log_2(K)$.
  \item $\Phi$ is $K$-to-1 {\rm ($\mu$-\ae)}.\hrulefill\ensuremath{\Box}
\end{list}
 \end{thm}

        \refstepcounter{thm}                     \begin{list}{} 			{\setlength{\leftmargin}{1em} 			\setlength{\rightmargin}{0em}}        \item {\bf Example \thethm:} Let $\lambda$ be the uniform Bernoulli measure on ${{\mathcal{ A}}^{\mathbb{Z}}}$.
Then $\lambda$ is invariant for any QGCA,
$\entr[\lambda]{\Phi}=\log_2(N)$, and $\Phi$ is $N$-to-1 ($\lambda$-\ae).
Indeed, $\lambda$ is the only $(\Phi, {{{\boldsymbol{\sigma}}}^{}} )$-invariant measure
with entropy $\log_2(N)$.  Thus, Proposition
\ref{Host.Maass.Martinez} is proved in \cite{HostMaassMartinez} by
first proving a special case of Theorem \ref{loop.ca.K.2.1} (when
$\Phi$ is a Ledrappier CA) and then showing that $K=N$.                   \hrulefill  \end{list}   			

Let $\mu$ be a measure on ${{\mathcal{ A}}^{\mathbb{Z}}}$.
If ${\mathfrak{ q}}$ is any partition of ${{\mathcal{ A}}^{\mathbb{Z}}}$, and ${\mathfrak{ S}}$ is any sigma-algebra,
define
 \begin{equation} 
\label{distropy}
\Entr{{\mathfrak{ q}}}{{\mathfrak{ S}}} \quad=\quad 
 \sum_{{\mathbf{ Q}}\in{\mathfrak{ q}}}   \int_{\mathbf{ Q}}
  \log_2\left(\Expct{{\mathbf{ Q}}}{{\mathfrak{ S}}}\rule[-0.5em]{0em}{1em}\right)({\mathbf{ x}}) \ d\mu[{\mathbf{ x}}].
 \end{equation} 
Let ${\mathfrak{ p}}_0$ be the partition of ${{\mathcal{ A}}^{\mathbb{Z}}}$ generated by zero-coordinate
cylinder sets,
and let ${\mathfrak{ p}}_{[\ell,n]} \ = \ \displaystyle \bigvee_{m=\ell}^n  {{{\boldsymbol{\sigma}}}^{-m}} ({\mathfrak{ p}}_0)$.  
Thus, ${\mathfrak{ B}}={\mathfrak{ p}}_{[-{\infty},{\infty}]}$ is the Borel sigma-algebra of ${{\mathcal{ A}}^{\mathbb{Z}}}$.
Let ${\mathfrak{ B}}^1 = \Phi^{-1}({\mathfrak{ B}})$.

  If $\mu$ is a $\Phi$-invariant measure, then $\displaystyle\entr{\Phi} \ =
\ \lim_{r{\rightarrow}{\infty}} \entr{\Phi,{\mathfrak{ p}}_{[-r,r]}}$, where
\begin{eqnarray*}
\entr{\Phi,{\mathfrak{ p}}_{[-r,r]}}
&=&
\Entr{{\mathfrak{ p}}_{[-r,r]}}{\bigvee_{t=1}^{\infty} \Phi^{-t}\left({\mathfrak{ p}}_{[-r,r]}\right)}.
\end{eqnarray*}

\begin{lemma}{\sf \label{entropy.lemma.1}}  
 If $\Phi:{{\mathcal{ A}}^{\mathbb{Z}}}{{\longrightarrow}}{{\mathcal{ A}}^{\mathbb{Z}}}$ is a QGCA, and $\mu$ is $(\Phi, {{{\boldsymbol{\sigma}}}^{}} )$-invariant, then  $\displaystyle \entr{\Phi}  =   \Entr{{\mathfrak{ p}}_0}{{\mathfrak{ B}}^1}$.
 \end{lemma}
\bprf
  Let ${\mathbf{ x}}\in{{\mathcal{ A}}^{\mathbb{Z}}}$ be an unknown sequence.  
  Because $\Phi$ is bipermutative, complete information about
$\left(\Phi^{t}({\mathbf{ x}})\right)_{\left[ -r,r \right]}$ (for $t\in{\left[ 0..T \right)}$)
is sufficient to reconstruct ${\mathbf{ x}}_{{\left[ -T-r,T+r \right]}}$, and vice versa.
In other words, we have an equality of partitions:
\[
 \bigvee_{t=0}^{T-1} \Phi^{-t}\left({\mathfrak{ p}}_{[-r,r]}\right)
\quad=\quad {\mathfrak{ p}}_{[-T-r,T+r]}.
\]
Letting $T{\rightarrow}{\infty}$, we get an equality of sigma-algebras:
$\displaystyle
 \bigvee_{t=0}^{\infty} \Phi^{-t}\left({\mathfrak{ p}}_{[-r,r]}\right)
\ = \ {\mathfrak{ p}}_{[-{\infty},{\infty}]} \ = \ {\mathfrak{ B}}$.
Applying $\Phi^{-1}$ to everything yields:
$\displaystyle \bigvee_{t=1}^{\infty} \Phi^{-t}\left({\mathfrak{ p}}_{[-r,r]}\right)
\ = \ \Phi^{-1}({\mathfrak{ B}}) \ = \ {\mathfrak{ B}}^1$. 
Hence,
\[
\entr{\Phi,{\mathfrak{ p}}_{[-r,r]}}
\quad=\quad
\Entr{{\mathfrak{ p}}_{[-r,r]}}{\bigvee_{t=1}^{\infty} \Phi^{-t}\left({\mathfrak{ p}}_{[-r,r]}\right)}
\quad=\quad
\Entr{{\mathfrak{ p}}_{[-r,r]}}{{\mathfrak{ B}}^{1}}
\]
 Now, $\Phi$ is bipermutative, so if we have complete
knowledge of $\Phi({\mathbf{ x}})$, then we can reconstruct ${\mathbf{ x}}$ from
knowledge only of $x_0$.  Thus,
$\Entr{{\mathfrak{ p}}_{[-r,r]}}{{\mathfrak{ B}}^{1}} \ = \ \Entr{{\mathfrak{ p}}_0}{{\mathfrak{ B}}^{1}}$.

Thus, \ 
$\entr{\Phi} \ = \ \displaystyle
\lim_{r{\rightarrow}{\infty}} \entr{\Phi,{\mathfrak{ p}}_{[-r,r]}}
\ = \ \lim_{r{\rightarrow}{\infty}}  \Entr{{\mathfrak{ p}}_0}{{\mathfrak{ B}}^{1}}
\ = \ \Entr{{\mathfrak{ p}}_0}{{\mathfrak{ B}}^{1}}$.
 {\tt \hrulefill $\Box$ } \end{list}  \medskip  

  For any ${\mathbf{ x}}\in{{\mathcal{ A}}^{\mathbb{Z}}}$, let ${\mathcal{ F}}\left({\mathbf{ x}}\right) \ = \ \Phi^{-1}\{\Phi({\mathbf{ x}})\}
\ = \ {\left\{ {\mathbf{ y}}\in{{\mathcal{ A}}^{\mathbb{Z}}} \; ; \; \Phi({\mathbf{ y}})=\Phi({\mathbf{ x}}) \right\} }$.  Hence, the sets ${\mathcal{ F}}\left({\mathbf{ x}}\right)$
(for ${\mathbf{ x}}\in{{\mathcal{ A}}^{\mathbb{Z}}}$) are the `minimal elements' of the sigma algebra
${\mathcal{ B}}^1$.

  The conditional expectation operator $\Expct{\bullet}{{\mathfrak{ B}}^1}$ defines
`fibre' measures $\mu_{\mathbf{ x}}$ (for ${{\forall}_{\mu}\;} {\mathbf{ x}}\in{{\mathcal{ A}}^{\mathbb{Z}}}$) having three properties:
\begin{description}
  \item[(F1)] For any measurable
${\mathbf{ U}}\subset{{\mathcal{ A}}^{\mathbb{Z}}}$ and for  ${{\forall}_{\mu}\;} {\mathbf{ x}}\in{{\mathcal{ A}}^{\mathbb{Z}}}$,
\quad $\mu_{\mathbf{ x}}({\mathbf{ U}}) \ = \  \Expct{{\mathbf{ U}}}{{\mathfrak{ B}}^1}({\mathbf{ x}})$.

  \item[(F2)] For any fixed ${\mathbf{ x}}\in{{\mathcal{ A}}^{\mathbb{Z}}}$, \quad $\mu_{\mathbf{ x}}$ is a 
probability measure
on ${{\mathcal{ A}}^{\mathbb{Z}}}$, and ${\sf supp}\left(\mu_{\mathbf{ x}}\right) \ = \ {\mathcal{ F}}\left({\mathbf{ x}}\right)$.

  \item[(F3)] For any fixed measurable ${\mathbf{ U}}\subset{{\mathcal{ A}}^{\mathbb{Z}}}$, the function
${{\mathcal{ A}}^{\mathbb{Z}}} \ni {\mathbf{ x}}\mapsto \mu_{\mathbf{ x}}({\mathbf{ U}}) \in {\mathbb{R}}$ is ${\mathfrak{ B}}^1$-measurable.
Hence, $\mu_{\mathbf{ x}} = \mu_{\mathbf{ y}}$ for any ${\mathbf{ y}}\in{\mathcal{ F}}\left({\mathbf{ x}}\right)$.
\end{description}

  Our goal is to show that there is some constant $K$ and, for
 ${{\forall}_{\mu}\;}{\mathbf{ x}}\in{{\mathcal{ A}}^{\mathbb{Z}}}$, there is a subset ${\mathcal{ E}}\subset{\mathcal{ F}}\left({\mathbf{ x}}\right)$ of
 cardinality $K$ so that $\mu_{\mathbf{ x}}$ is uniformly distributed on ${\mathcal{ E}}$.

\begin{lemma}{\sf \label{mu.shift.cov}}  For any measurable
${\mathbf{ U}}\subset{{\mathcal{ A}}^{\mathbb{Z}}}$ and for  ${{\forall}_{\mu}\;} {\mathbf{ x}}\in{{\mathcal{ A}}^{\mathbb{Z}}}$,\quad
 $\mu_{\mathbf{ x}}\left( {{{\boldsymbol{\sigma}}}^{-1}} ({\mathbf{ U}})\right) \ = \ \mu_{ {{{\boldsymbol{\sigma}}}^{}} ({\mathbf{ x}})}({\mathbf{ U}})$.
 \end{lemma}
\bprf
For ${{\forall}_{\mu}\;} {\mathbf{ x}}\in{{\mathcal{ A}}^{\mathbb{Z}}}$, property {\bf(F1)} says
\begin{eqnarray*}
\mu_{\mathbf{ x}}\left( {{{\boldsymbol{\sigma}}}^{-1}} ({\mathbf{ U}})\right) & = & 
                 \Expct{ {{{\boldsymbol{\sigma}}}^{-1}} ({\mathbf{ U}})}{{\mathfrak{ B}}^1}({\mathbf{ x}}),\\
\mbox{\ and \ } \quad \mu_{ {{{\boldsymbol{\sigma}}}^{}} ({\mathbf{ x}})}({\mathbf{ U}})
                  & = & \Expct{{\mathbf{ U}}}{{\mathfrak{ B}}^1}\left( {{{\boldsymbol{\sigma}}}^{}} ({\mathbf{ x}})\right).
\end{eqnarray*}
   Thus, we must show that
$\Expct{ {{{\boldsymbol{\sigma}}}^{-1}} ({\mathbf{ U}})}{{\mathfrak{ B}}^1}({\mathbf{ x}}) \ = \
 \Expct{{\mathbf{ U}}}{{\mathfrak{ B}}^1}\left( {{{\boldsymbol{\sigma}}}^{}} ({\mathbf{ x}})\right)$ for ${{\forall}_{\mu}\;} {\mathbf{ x}}\in{{\mathcal{ A}}^{\mathbb{Z}}}$.
But $\Expct{ {{{\boldsymbol{\sigma}}}^{-1}} ({\mathbf{ U}})}{{\mathfrak{ B}}^1}$ and $\Expct{{\mathbf{ U}}}{{\mathfrak{ B}}^1}$ are
${\mathfrak{ B}}^1$-measurable functions, so it suffices to show that
$\displaystyle \int_{\mathbf{ B}} \Expct{ {{{\boldsymbol{\sigma}}}^{-1}} ({\mathbf{ U}})}{{\mathfrak{ B}}^1}({\mathbf{ X}}) \ d\mu[{\mathbf{ x}}]
\ = \ \int_{\mathbf{ B}} \Expct{{\mathbf{ U}}}{{\mathfrak{ B}}^1}\left( {{{\boldsymbol{\sigma}}}^{}} ({\mathbf{ x}})\right) \ d\mu[{\mathbf{ x}}]$,
 for any ${\mathbf{ B}}\in{\mathfrak{ B}}^1$. But
\begin{eqnarray*}
 \int_{\mathbf{ B}} \Expct{ {{{\boldsymbol{\sigma}}}^{-1}} ({\mathbf{ U}})}{{\mathfrak{ B}}^1}({\mathbf{ x}}) \ d\mu[{\mathbf{ x}}]
&\displaystyle\raisebox{-1ex}{$\overline{\overline{{\scriptscriptstyle{\mathrm{(E)}}}}}$}& 
 \int_{\mathbf{ B}} {{{\mathsf{ 1\!\!1}}}_{{ {{{\boldsymbol{\sigma}}}^{-1}} ({\mathbf{ U}})}}}({\mathbf{ x}}) \ d\mu[{\mathbf{ x}}]
\quad=\quad
  \mu\left[{\mathbf{ B}}\cap  {{{\boldsymbol{\sigma}}}^{-1}} ({\mathbf{ U}})\right]
\\&\displaystyle\raisebox{-1ex}{$\overline{\overline{{\scriptscriptstyle{\mathrm{(I)}}}}}$}&
  \mu\left[ {{{\boldsymbol{\sigma}}}^{}} ({\mathbf{ B}})\cap {\mathbf{ U}}\right]
\quad=\quad
 \int_{ {{{\boldsymbol{\sigma}}}^{}} ({\mathbf{ B}})} {{{\mathsf{ 1\!\!1}}}_{{{\mathbf{ U}}}}}({\mathbf{ x}}') \ d\mu[{\mathbf{ x}}']
\\&\displaystyle\raisebox{-1ex}{$\overline{\overline{{\scriptscriptstyle{\mathrm{(E)}}}}}$}&
\int_{ {{{\boldsymbol{\sigma}}}^{}} ({\mathbf{ B}})}  \Expct{{\mathbf{ U}}}{{\mathfrak{ B}}^1}({\mathbf{ x}}') \ d\mu[{\mathbf{ x}}']
\\&\displaystyle\raisebox{-1ex}{$\overline{\overline{{\scriptscriptstyle{\mathrm{(S)}}}}}$}&
\int_{{\mathbf{ B}}}  \Expct{{\mathbf{ U}}}{{\mathfrak{ B}}^1}\left( {{{\boldsymbol{\sigma}}}^{}} ({\mathbf{ x}})\right) \ d\mu[{\mathbf{ x}}],
\end{eqnarray*}
as desired.  Here {\bf(E)} is the defining property of conditional
expectation, {\bf(I)} is because $\mu$ is $ {{{\boldsymbol{\sigma}}}^{}} $-invariant, and {\bf(S)}
is the substitution ${\mathbf{ x}}'= {{{\boldsymbol{\sigma}}}^{}} ({\mathbf{ x}})$ (again because $\mu$ is
$ {{{\boldsymbol{\sigma}}}^{}} $-invariant).
 {\tt \hrulefill $\Box$ } \end{list}  \medskip  

For any ${\mathbf{ x}}\in{{\mathcal{ A}}^{\mathbb{Z}}}$, let $\eta({\mathbf{ x}}) = \mu_{\mathbf{ x}}\{{\mathbf{ x}}\}$.  Thus, if
${\mathbf{ y}}$ is an unknown, $\mu$-random sequence, then $\eta({\mathbf{ x}})$ represents the
conditional probability that ${\mathbf{ y}}={\mathbf{ x}}$, given that
$\Phi({\mathbf{ y}})=\Phi({\mathbf{ x}})$.

\begin{lemma}{\sf \label{eta.shift.inv}}  
 \setcounter{enumi}{\thethm} \begin{list}{{\bf (\alph{enumii})}}{\usecounter{enumii}} 			{\setlength{\leftmargin}{0em} 			\setlength{\rightmargin}{0em}}
  \item  $\eta$ is  $ {{{\boldsymbol{\sigma}}}^{}} $-invariant {\rm ($\mu$-\ae)}.

  \item If $\mu$ is $ {{{\boldsymbol{\sigma}}}^{}} $-ergodic, then  
there is some $H\in{\mathbb{R}}$ such that $\eta({\mathbf{ x}})=H$ for ${{\forall}_{\mu}\;} {\mathbf{ x}}\in{{\mathcal{ A}}^{\mathbb{Z}}}$.
 
  \item If $\mu$ is also $\Phi$-invariant, then $\eta$ is $\Phi$-invariant {\rm ($\mu$-\ae)}.
\end{list}
 \end{lemma}
\bprf 
{\bf(a)} \quad  $\eta\left( {{{\boldsymbol{\sigma}}}^{}} ({\mathbf{ x}})\right) \ = \ \mu_{ {{{\boldsymbol{\sigma}}}^{}} ({\mathbf{ x}})} \{ {{{\boldsymbol{\sigma}}}^{}} ({\mathbf{ x}})\}
\ \ \displaystyle\raisebox{-1ex}{$\overline{\overline{{\scriptscriptstyle{\mathrm{(L\ref{mu.shift.cov})}}}}}$} \ \ 
 \mu_{{\mathbf{ x}}} \left( {{{\boldsymbol{\sigma}}}^{-1}} \{ {{{\boldsymbol{\sigma}}}^{}} ({\mathbf{ x}})\}\right)
\ \ \displaystyle\raisebox{-1ex}{$\overline{\overline{{\scriptscriptstyle{\mathrm{(\dagger)}}}}}$} \ \ \mu_{{\mathbf{ x}}} \{{\mathbf{ x}}\}
\ = \ \eta({\mathbf{ x}})$. 

 {\bf(L\ref{mu.shift.cov})} is Lemma \ref{mu.shift.cov}.\quad
  $(\dagger)$ is because $ {{{\boldsymbol{\sigma}}}^{}} $ is invertible on ${{\mathcal{ A}}^{\mathbb{Z}}}$.
Parts {\bf(b)}  and {\bf(c)}  follow.  {\tt \hrulefill $\Box$ } \end{list}  \medskip  

\begin{lemma}{\sf \label{entropy.lemma.2}}  
  If $\mu$ is $ {{{\boldsymbol{\sigma}}}^{}} $-ergodic, then $\entr{\Phi} \ = \  -\log_2(H)$.
 \end{lemma}
\bprf
Lemma \ref{entropy.lemma.1} and eqn. (\ref{distropy}) imply:
$\displaystyle \entr{\Phi} \ = \ 
 - \sum_{{\mathbf{ P}}\in{\mathfrak{ p}}_0}   \int_{\mathbf{ P}}
  \log_2\left(\Expct{{\mathbf{ P}}}{{\mathfrak{ B}}^1}\rule[-0.5em]{0em}{1em}\right)({\mathbf{ x}}) \ d\mu[{\mathbf{ x}}]$.

\refstepcounter{claimcount}                {\bf Claim \theclaimcount: \ }{\sl  \label{entropy.lemma.2.1}
For all ${\mathbf{ P}}\in{\mathfrak{ p}}_0$, and for ${{\forall}_{\mu}\;} {\mathbf{ x}}\in{\mathbf{ P}}$,
\quad $ \Expct{{\mathbf{ P}}}{{\mathfrak{ B}}^1}({\mathbf{ x}}) \ = \  H$.}
\bprf
$ \Expct{{\mathbf{ P}}}{{\mathfrak{ B}}^1}({\mathbf{ x}})
\ \  \displaystyle\raisebox{-1ex}{$\overline{\overline{{\scriptscriptstyle{\mathrm{(F1)}}}}}$} \ \ 
\mu_{\mathbf{ x}}({\mathbf{ P}}) 
\ \  \displaystyle\raisebox{-1ex}{$\overline{\overline{{\scriptscriptstyle{\mathrm{(F2)}}}}}$} \ \
\mu_{\mathbf{ x}}\left({\mathbf{ P}}\cap{\mathcal{ F}}\left({\mathbf{ x}}\right)\rule[-0.5em]{0em}{1em}\right) 
\ \ \displaystyle\raisebox{-1ex}{$\overline{\overline{{\scriptscriptstyle{\mathrm{(c2)}}}}}$} \ \
\mu_{\mathbf{ x}}\{{\mathbf{ x}}\}
\ = \ \eta({\mathbf{ x}})
\ \ \displaystyle\raisebox{-1ex}{$\overline{\overline{{\scriptscriptstyle{\mathrm{(\ref{eta.shift.inv}b)}}}}}$} \ \ H$.

Here, {\bf(c2)} follows from Claim 2 below,
and {\bf(\ref{eta.shift.inv}b)} is by Corollary \ref{eta.shift.inv}{\bf(b)}.
 {\tt \dotfill~$\Diamond$~[Claim~\theclaimcount] }\end{list}

\refstepcounter{claimcount}                {\bf Claim \theclaimcount: \ }{\sl  \label{entropy.lemma.2.2}
If ${\mathbf{ x}}\in{\mathbf{ P}}\in{\mathfrak{ p}}_0$, then  ${\mathbf{ P}}\cap{\mathcal{ F}}\left({\mathbf{ x}}\right) = \{{\mathbf{ x}}\}$.}
\bprf
 $\Phi$ is bipermutative, so
if ${\mathbf{ y}}\in{\mathcal{ F}}\left({\mathbf{ x}}\right)$, then ${\mathbf{ y}}$ is entirely determined by $y_0$.
Thus, $\left( \ \rule[-0.5em]{0em}{1em}       \begin{minipage}{40em}       \begin{tabbing}         ${\mathbf{ y}}\in{\mathbf{ P}}$        \end{tabbing}      \end{minipage} \ \right)\ensuremath{\Longrightarrow}\left( \ \rule[-0.5em]{0em}{1em}       \begin{minipage}{40em}       \begin{tabbing}         $y_0=x_0$        \end{tabbing}      \end{minipage} \ \right)
\ensuremath{\Longrightarrow}\left( \ \rule[-0.5em]{0em}{1em}       \begin{minipage}{40em}       \begin{tabbing}         ${\mathbf{ y}}={\mathbf{ x}}$        \end{tabbing}      \end{minipage} \ \right)$.
\hrulefill $\Diamond$ {\tt[Claim~\theclaimcount]~~$\Box$} \end{list} \end{list}  \medskip  

  Our goal is to show that $H=\frac{1}{K}$
for some $K$.

  Suppose $|{\mathcal{ A}}| = N$, and identify ${\mathcal{ A}}$ with the group ${{\mathbb{Z}}_{/N}}$
in an arbitrary way.  Define $\tau:{{\mathcal{ A}}^{\mathbb{Z}}}{{\longrightarrow}}{{\mathcal{ A}}^{\mathbb{Z}}}$ as follows.
For any ${\mathbf{ x}}\in{{\mathcal{ A}}^{\mathbb{Z}}}$, $\tau({\mathbf{ x}})={\mathbf{ y}}$, where ${\mathbf{ y}}$ is the unique
element in ${\mathcal{ F}}\left({\mathbf{ x}}\right)$ such that $y_0 = x_0 + 1 \pmod{N}$. 
Existence/uniqueness of ${\mathbf{ y}}$ follows from  bipermutativity.

Note that $\tau(\mu)\neq \mu$, so a statement which is true
$\mu$-\ae\ may {\em not} be true $\tau(\mu)$-\ae.
For example, Lemma \ref{eta.shift.inv}{\bf(c)} does
{\em not} imply that $\eta\left(\Phi\left[\tau({\mathbf{ x}})\rule[-0.5em]{0em}{1em}\right]\right) \ = \
\eta\left(\tau({\mathbf{ x}})\rule[-0.5em]{0em}{1em}\right)$ for ${{\forall}_{\mu}\;} {\mathbf{ x}}$.

 Let ${\mathbf{ E}}_n \ = \ {\left\{ {\mathbf{ x}}\in{{\mathcal{ A}}^{\mathbb{Z}}} \; ; \; \eta\left(\tau^n({\mathbf{ x}})\rule[-0.5em]{0em}{1em}\right) > \ 0 \right\} }$.
Let $\mu_n \ = \ \tau^n\left( {{{\mathsf{ 1\!\!1}}}_{{{\mathbf{ E}}_n}}} \cdot \mu\right)$.

\begin{lemma}{\sf \label{mu.n.ll.mu}}  
  $\mu_n$ is absolutely continuous relative to $\mu$.
 \end{lemma}
\bprf
  Suppose ${\mathbf{ Z}}\subset{{\mathcal{ A}}^{\mathbb{Z}}}$ is Borel-measurable, and $\mu[{\mathbf{ Z}}]=0$.
We want to show $\mu_n[{\mathbf{ Z}}]=0$ also.  But
\  $\mu_n[{\mathbf{ Z}}] \ = \  \left({{{\mathsf{ 1\!\!1}}}_{{{\mathbf{ E}}_n}}} \cdot \mu\right)\left[ \tau^{-n}({\mathbf{ Z}})\right]
 \ = \  \mu\left[\tau^{-n}({\mathbf{ Z}}) \cap {\mathbf{ E}}_n\right]$, \ 
so it suffices to show:
\refstepcounter{claimcount}                {\bf Claim \theclaimcount: \ }{\sl  For ${{\forall}_{\mu}\;} {\mathbf{ z}}\in\tau^{-n}({\mathbf{ Z}})$, \quad $\eta\left(\tau^n({\mathbf{ z}})\right) \ = 0$;
\ hence ${\mathbf{ z}}\not\in{\mathbf{ E}}_n$.}
\bprf
$\displaystyle \int_{{\mathcal{ A}}^{\mathbb{Z}}} \mu_{\mathbf{ x}}[{\mathbf{ Z}}] \ { \;\; d\mu}[{\mathbf{ x}}] \ = \  \mu[{\mathbf{ Z}}] \ = \ 0$.  Hence,
for ${{\forall}_{\mu}\;} {\mathbf{ x}}\in{{\mathcal{ A}}^{\mathbb{Z}}}$, 
\quad
$  \mu_{\mathbf{ x}}[{\mathbf{ Z}}]  \ = \ 0$. \
But if ${\mathbf{ z}}\in\tau^{-n}({\mathbf{ Z}})$, then 
$\tau^n({\mathbf{ z}})\in{\mathcal{ F}}({\mathbf{ z}})\cap{\mathbf{ Z}}$, so we get 
$\eta\left(\tau^n({\mathbf{ z}})\right) \ = \ \mu_{\tau^n({\mathbf{ z}})}\{\tau^n({\mathbf{ z}})\} \ = \
\mu_{\mathbf{ z}}\{\tau^n({\mathbf{ z}})\} \ \leq \  \mu_{\mathbf{ z}}[{\mathbf{ Z}}\cap {\mathcal{ F}}({\mathbf{ x}})] \ \leq \ \mu_{\mathbf{ z}}[{\mathbf{ Z}}] \ = \ 0$. 
\hrulefill $\Diamond$ {\tt[Claim~\theclaimcount]~~$\Box$} \end{list} \end{list}  \medskip

\begin{cor}{\sf \label{eta.constant.2}}  
 If $\mu$ is $ {{{\boldsymbol{\sigma}}}^{}} $-ergodic and  $\Phi$-invariant, 
then $\eta$ is $\Phi$-invariant {\rm($\mu_n$-\ae)}.
 \end{cor}
\bprf
Lemma \ref{mu.n.ll.mu} means that
 a statement which is true for ${{\forall}_{\mu}\;}
{\mathbf{ x}}$ is also true for $\forall_{\mu_n} {\mathbf{ x}}$. 
Now apply Lemma \ref{eta.shift.inv}{\bf(c)}.
 {\tt \hrulefill $\Box$ } \end{list}  \medskip  

\begin{cor}{\sf \label{eta.constant.3}}  
  For ${{\forall}_{\mu}\;} {\mathbf{ x}}\in {\mathbf{ E}}_n$, \quad $\eta({\mathbf{ x}})
\ = \ \eta\left(\tau^n({\mathbf{ x}})\rule[-0.5em]{0em}{1em}\right)$.
 \end{cor}
\bprf
$\eta\left({\mathbf{ x}}\right)
\quad\displaystyle\raisebox{-1ex}{$\overline{\overline{{\scriptscriptstyle{\mathrm{(\ref{eta.shift.inv}c)}}}}}$} \quad \eta\left(\Phi[{\mathbf{ x}}]\rule[-0.5em]{0em}{1em}\right)
 \quad \displaystyle\raisebox{-1ex}{$\overline{\overline{{\scriptscriptstyle{\mathrm{(*)}}}}}$} \quad \eta\left(\Phi\left[\tau^n({\mathbf{ x}})\rule[-0.5em]{0em}{1em}\right]\right)
 \quad \displaystyle\raisebox{-1ex}{$\overline{\overline{{\scriptscriptstyle{\mathrm{(C\ref{eta.constant.2})}}}}}$} \quad \eta\left(\tau^n({\mathbf{ x}})\rule[-0.5em]{0em}{1em}\right)$.
Here, {\bf(\ref{eta.shift.inv}c)} is by Corollary \ref{eta.shift.inv}{\bf(c)},
$(*)$ is because $\tau^n({\mathbf{ x}})\in{\mathcal{ F}}\left({\mathbf{ x}}\right)$, and
{\bf(C\ref{eta.constant.2})} is by Corollary \ref{eta.constant.2}.
 {\tt \hrulefill $\Box$ } \end{list}  \medskip  

  Now, let ${\mathcal{ E}}({\mathbf{ x}}) = {\left\{ {\mathbf{ y}}\in{\mathcal{ F}}\left({\mathbf{ x}}\right) \; ; \; \eta({\mathbf{ y}})>0 \right\} }$.  

\begin{cor}{\sf \label{eta.constant.4}}  
   For ${{\forall}_{\mu}\;} {\mathbf{ x}}\in{{\mathcal{ A}}^{\mathbb{Z}}}$, \quad $\mu_{\mathbf{ x}}$ is equidistributed
on ${\mathcal{ E}}({\mathbf{ x}})$.  If ${\sf card}\left[{\mathcal{ E}}({\mathbf{ x}})\right]=K$, then $\mu_{\mathbf{ x}}$ assigns
mass $\frac{1}{K}$ to each element in ${\mathcal{ E}}({\mathbf{ x}})$.  In particular,
$\eta({\mathbf{ x}}) \ = \ \frac{1}{K}$.
 \end{cor}
\bprf
  By definition, \
$\displaystyle  1 \ \  \displaystyle\raisebox{-1ex}{$\overline{\overline{{\scriptscriptstyle{\mathrm{(F2)}}}}}$} \ \  \mu_{\mathbf{ x}} \left(\rule[-0.5em]{0em}{1em}{\mathcal{ F}}\left({\mathbf{ x}}\right)\right)
      = \  \sum_{{\mathbf{ y}}\in{\mathcal{ F}}\left({\mathbf{ x}}\right)} \mu_{\mathbf{ x}} \{{\mathbf{ y}}\}
    \ \  \displaystyle\raisebox{-1ex}{$\overline{\overline{{\scriptscriptstyle{\mathrm{(F3)}}}}}$} \ \  \sum_{{\mathbf{ y}}\in{\mathcal{ F}}\left({\mathbf{ x}}\right)} \mu_{\mathbf{ y}} \{{\mathbf{ y}}\}
      = \  \sum_{{\mathbf{ y}}\in{\mathcal{ E}}({\mathbf{ x}})} \eta({\mathbf{ y}})$.

  However, if ${\mathbf{ y}}=\tau^n({\mathbf{ x}})$, then
$\left( \ \rule[-0.5em]{0em}{1em}       \begin{minipage}{40em}       \begin{tabbing}         ${\mathbf{ y}}\in{\mathcal{ E}}({\mathbf{ x}})$        \end{tabbing}      \end{minipage} \ \right) \iff \left( \ \rule[-0.5em]{0em}{1em}       \begin{minipage}{40em}       \begin{tabbing}         ${\mathbf{ x}}\in{\mathbf{ E}}_n$        \end{tabbing}      \end{minipage} \ \right)$,
in which case Corollary \ref{eta.constant.3} implies that
$\eta({\mathbf{ y}}) = \eta({\mathbf{ x}})$.  Hence,  \ 
$\displaystyle  1 \  = \ \sum_{{\mathbf{ y}}\in{\mathcal{ E}}({\mathbf{ x}})} \eta({\mathbf{ y}}) 
\ = \  \sum_{{\mathbf{ y}}\in{\mathcal{ E}}({\mathbf{ x}})} \eta({\mathbf{ x}})
 \ = \  K\cdot \eta({\mathbf{ x}})$, \
  where $K={\sf card}\left[{\mathcal{ E}}({\mathbf{ x}})\right]$.  We conclude that $\eta({\mathbf{ x}})=\frac{1}{K}$.
 {\tt \hrulefill $\Box$ } \end{list}  \medskip  

\begin{cor}{\sf \label{eta.constant.5}}  
  There is some $K$ so that, for ${{\forall}_{\mu}\;} {\mathbf{ x}}\in{{\mathcal{ A}}^{\mathbb{Z}}}$, \quad 
${\sf card}\left[{\mathcal{ E}}({\mathbf{ x}})\right]=K$, and $\mu_{\mathbf{ x}}$ assigns
mass $\frac{1}{K}$ to each element of ${\mathcal{ E}}({\mathbf{ x}})$.
Thus,  $H = \frac{1}{K}$. Thus, $\entr{\Phi} = \log_2(K)$.
 \end{cor}
\bprf Combine Corollaries \ref{eta.shift.inv}{\bf(b)}
 and \ref{eta.constant.4}.
Then apply Lemma \ref{entropy.lemma.2}.
 {\tt \hrulefill $\Box$ } \end{list}  \medskip

\bprf[Proof of Theorem \ref{loop.ca.K.2.1}:]
  Let ${\mathcal{ U}} = {\left\{ {\mathbf{ x}}\in{{\mathcal{ A}}^{\mathbb{Z}}} \; ; \; {\sf card}\left[{\mathcal{ E}}({\mathbf{ x}})\right]=K \right\} }$.  Then
Corollary \ref{eta.constant.5} says $\mu({\mathcal{ U}})=1$. Since
$\mu$ is $\Phi$-invariant, it follows that $\mu(\Phi^{-1}({\mathcal{ U}}))=1$
also; hence $\Phi^{-1}({\mathcal{ U}}) \ \displaystyle\raisebox{-0.6ex}{$\overline{\overline{{\scriptstyle{\mathrm{\,\mu\,}}}}}$} \ {\mathcal{ U}}$.

  Thus, for ${{\forall}_{\mu}\;} {\mathbf{ u}}\in{\mathcal{ U}}$, there is some ${\mathbf{ x}}\in{\mathcal{ U}}$ so that
$\Phi({\mathbf{ x}})={\mathbf{ u}}$.  But then $\Phi^{-1}({\mathbf{ u}})={\mathcal{ F}}\left({\mathbf{ x}}\right)$, and
$\Phi^{-1}({\mathbf{ u}})\cap{\mathcal{ U}} = {\mathcal{ F}}\left({\mathbf{ x}}\right)\cap{\mathcal{ U}} \ = \
{\mathcal{ E}}({\mathbf{ x}})$ is a set of cardinality $K$,
by definition of ${\mathcal{ U}}$.
 {\tt \hrulefill $\Box$ } \end{list}  \medskip

\section{\label{S:eca} Endomorphic Cellular Automata}

  A {\bf group shift} is a sequence space ${\mathcal{ A}}^{\mathbb{Z}}$ equipped with a
topological group structure such that $ {{{\boldsymbol{\sigma}}}^{}} $ is a group automorphism.
Equivalently, the multiplication operation $\bullet$ on ${\mathcal{ A}}^{\mathbb{Z}}$ is
defined by some {\bf local multiplication map} $\psi:{\mathcal{ A}}^{\left[ -\ell..r \right]}
\times {\mathcal{ A}}^{\left[ -\ell..r \right]}{{\longrightarrow}}{\mathcal{ A}}$ so that, if ${\mathbf{ a}},{\mathbf{ b}}\in{\mathcal{ A}}^{\mathbb{Z}}$
and ${\mathbf{ c}} = {\mathbf{ a}}\bullet {\mathbf{ b}}$, then $c_0 \ = \ \psi(a_{-\ell},\ldots,a_r; \ 
b_{-\ell},\ldots,b_{r})$.

  The most obvious group shift is a {\em product group}, where ${\mathcal{ A}}$ is
a finite group and  multiplication on ${\mathcal{ A}}^{\mathbb{Z}}$ is defined
componentwise.  However, this is not the only group shift \cite{Kitchens}.

  An {\bf endomorphic cellular automaton} (ECA) is a cellular automaton
$\Phi:{\mathcal{ A}}^{\mathbb{Z}}{{\longrightarrow}}{\mathcal{ A}}^{\mathbb{Z}}$ which is also a group endomorphism of
${\mathcal{ A}}^{\mathbb{Z}}$.  For example, it is easy to verify:

\begin{prop}{\sf \label{biperm.eca.abel.prod.grp}}  
  Let $({\mathcal{ A}},+)$ be an additive abelian group.  Let ${\mathcal{ A}}^{\mathbb{Z}}$ be
the product group.  Let $\Phi:{\mathcal{ A}}^{\mathbb{Z}}{{\longrightarrow}}{\mathcal{ A}}^{\mathbb{Z}}$ be a
right-sided, nearest neighbour CA, with local map
$\phi:{\mathcal{ A}}^{\{0,1\}}{{\longrightarrow}}{\mathcal{ A}}$.  Then:
 \setcounter{enumi}{\thethm} \begin{list}{{\bf (\alph{enumii})}}{\usecounter{enumii}} 			{\setlength{\leftmargin}{0em} 			\setlength{\rightmargin}{0em}}
  \item $\Phi$ is an ECA iff $\phi(a_0,a_1) = \phi_0(a_0) + \phi_1(a_1)$,
where $\phi_0, \phi_1$ are endomorphisms of ${\mathcal{ A}}$.

  \item $\Phi$ is bipermutative iff $\phi_0$ and $\phi_1$ are automorphisms
of ${\mathcal{ A}}$.
\hrulefill\ensuremath{\Box}
\end{list}
 \end{prop}

  We will now apply the results of 
\S\ref{S:degree} to bipermutative ECA, to prove:

\begin{thm}{\sf \label{biperm.eca.measure}}  
  Let ${\mathcal{ A}}^{\mathbb{Z}}$ be a group shift and let $\Phi:{\mathcal{ A}}^{\mathbb{Z}}{{\longrightarrow}}{\mathcal{ A}}^{\mathbb{Z}}$ be
a bipermutative ECA. Suppose $\ker(\Phi)$
 contains no nontrivial $ {{{\boldsymbol{\sigma}}}^{}} $-invariant subgroups.

If $\mu$ is $\Phi$-invariant and totally $ {{{\boldsymbol{\sigma}}}^{}} $-ergodic, and
$\entr{\Phi}>0$, then $\mu=\lambda$. \hrulefill\ensuremath{\Box}  \end{thm}

Recall from \S\ref{S:degree} that if ${\mathbf{ x}}\in{\mathcal{ A}}^{\mathbb{Z}}$, then ${\mathcal{ F}}\left({\mathbf{ x}}\right) = \Phi^{-1}\{\Phi({\mathbf{ x}})\}$.

\begin{lemma}{\sf \label{kernel.biperm.eca.1}}  
  Let $\Phi:{\mathcal{ A}}^{\mathbb{Z}}{{\longrightarrow}}{\mathcal{ A}}^{\mathbb{Z}}$ be
a bipermutative ECA on a group shift.  Let ${\mathcal{ K}}  =\ker(\Phi)$.
 \setcounter{enumi}{\thethm} \begin{list}{{\bf (\alph{enumii})}}{\usecounter{enumii}} 			{\setlength{\leftmargin}{0em} 			\setlength{\rightmargin}{0em}}
  \item For any ${\mathbf{ x}}\in{\mathcal{ A}}^{\mathbb{Z}}$,\quad ${\mathcal{ F}}\left({\mathbf{ x}}\right) = {\mathbf{ x}} \bullet{\mathcal{ K}}$.

  \item  Let ${\mathbf{ e}}\in{\mathcal{ A}}^{\mathbb{Z}}$ be the identity element.  Then ${\mathbf{ e}}$ is
a constant sequence ---ie. there is some $e\in{\mathcal{ A}}$ so that 
${\mathbf{ e}}=(....,e,e,e,....)$.

  \item ${\mathcal{ K}}$ is $ {{{\boldsymbol{\sigma}}}^{}} $-invariant.  Also, if ${\mathbf{ k}}\in{\mathcal{ K}}$, then
  ${\mathbf{ k}}$ is entirely determined by $k_0$.

  \item There is a natural bijection $\zeta:{\mathcal{ A}}{{\longrightarrow}}{\mathcal{ K}}$, where $\zeta[a]$
is the unique element ${\mathbf{ k}}\in{\mathcal{ K}}$ with $k_0=a$.   In particular, 
$\zeta[e]={\mathbf{ e}}$.

  \item  There is a permutation
$\rho:{\mathcal{ A}}{{\longrightarrow}}{\mathcal{ A}}$ so that $ {{{\boldsymbol{\sigma}}}^{}} \left(\zeta[a]\right) \ = \ \zeta\left[
\rho(a)\right]$.  In particular, $\rho(e)=e$. 

  It follows that every element of ${\mathcal{ K}}$ is $P$-periodic, for some $P<|{\mathcal{ A}}|$.

  \item Any $ {{{\boldsymbol{\sigma}}}^{}} $-invariant subgroup ${\mathcal{ J}}\prec{\mathcal{ K}}$ is thus a disjoint
union of periodic $ {{{\boldsymbol{\sigma}}}^{}} $-orbits, which corresponds to a disjoint union of
$\rho$-orbits in ${\mathcal{ A}}$. 

  \item In particular: 
\quad $\left( \ \rule[-0.5em]{0em}{1em}       \begin{minipage}{40em}       \begin{tabbing}         ${\mathcal{ A}}\setminus\{e\}$ consists of a single
$\rho$-orbit        \end{tabbing}      \end{minipage} \ \right)$

$\iff
\left( \ \rule[-0.5em]{0em}{1em}       \begin{minipage}{40em}       \begin{tabbing}         ${\mathcal{ K}}$ has no nontrivial $ {{{\boldsymbol{\sigma}}}^{}} $-invariant subgroups        \end{tabbing}      \end{minipage} \ \right).$
\end{list}
 \end{lemma}
\bprf
  {\bf(a)} is a basic property of group homomorphisms.  To see
{\bf(b)}, recall that $ {{{\boldsymbol{\sigma}}}^{}} $ is a group automorphism of ${\mathcal{ A}}^{\mathbb{Z}}$.
Thus, $ {{{\boldsymbol{\sigma}}}^{}} ({\mathbf{ e}})={\mathbf{ e}}$, so ${\mathbf{ e}}$ must be constant.
{\bf(c)} follows from {\bf(b)} and the fact that $\Phi$ is bipermutative.
Then {\bf(c)} implies {\bf(d)} implies {\bf(e)} implies {\bf(f)}.
 {\tt \hrulefill $\Box$ } \end{list}  \medskip  

  If $({\mathcal{ A}},+)$ is abelian and ${\mathcal{ A}}^{\mathbb{Z}}$ is the product group, then
Lemma \ref{kernel.biperm.eca.1} takes the form:

\begin{lemma}{\sf \label{kernel.biperm.eca.2}}  
  Let $({\mathcal{ A}},+)$ be an abelian group and let ${\mathcal{ A}}^{\mathbb{Z}}$ 
be the product group. Let $\Phi:{\mathcal{ A}}^{\mathbb{Z}}{{\longrightarrow}}{\mathcal{ A}}^{\mathbb{Z}}$ be
a bipermutative ECA and let ${\mathcal{ K}} =\ker(\Phi)$. 
 \setcounter{enumi}{\thethm} \begin{list}{{\bf (\alph{enumii})}}{\usecounter{enumii}} 			{\setlength{\leftmargin}{0em} 			\setlength{\rightmargin}{0em}}
  \item For any ${\mathbf{ x}}\in{\mathcal{ A}}^{\mathbb{Z}}$,\quad ${\mathcal{ F}}\left({\mathbf{ x}}\right) = {\mathbf{ x}} + {\mathcal{ K}}$.

  \item The map $\zeta:{\mathcal{ A}}{{\longrightarrow}}{\mathcal{ K}}$ from Lemma
 \ref{kernel.biperm.eca.1}{\bf(d)} is a group isomorphism.

  \item  The map $\rho:{\mathcal{ A}}{{\longrightarrow}}{\mathcal{ A}}$ from Lemma 
\ref{kernel.biperm.eca.1}{\bf(e)} is a group automorphism.  To be precise,
suppose $\Phi$ has local map $\phi(a_0,a_1) = \phi_0(a_0) + \phi_1(a_1)$,
where $\phi_0$ and $\phi_1$ are automorphisms of ${\mathcal{ A}}$, as in
Proposition \ref{biperm.eca.abel.prod.grp}{\bf(b)}.  Then $\rho
\ = \ -\phi_1^{-1}\circ \phi_0$.

  \item If ${\mathcal{ J}}\prec{\mathcal{ K}}$ is a $ {{{\boldsymbol{\sigma}}}^{}} $-invariant subgroup, then
${\mathcal{ J}} = \zeta({\mathcal{ B}})$, where ${\mathcal{ B}}\prec{\mathcal{ A}}$ is a $\rho$-invariant 
subgroup of ${\mathcal{ A}}$.

  \item In particular,
 \quad $\left( \ \rule[-0.5em]{0em}{1em}       \begin{minipage}{40em}       \begin{tabbing}         ${\mathcal{ A}}$ has no nontrivial $\rho$-invariant
subgroups        \end{tabbing}      \end{minipage} \ \right)$

$\iff
\left( \ \rule[-0.5em]{0em}{1em}       \begin{minipage}{40em}       \begin{tabbing}         ${\mathcal{ K}}$ has no nontrivial $ {{{\boldsymbol{\sigma}}}^{}} $-invariant subgroups        \end{tabbing}      \end{minipage} \ \right).$
\end{list}
 \end{lemma}
\bprf
  We need only verify the claim in {\bf(b)} that $\zeta$ is a group
homomorphism.  To see this,  suppose ${\mathbf{ k}} = \zeta(a)$ and ${\mathbf{ k}}' = \zeta(a')$.
Let ${\mathbf{ j}} = {\mathbf{ k}}+{\mathbf{ k}}'$ and let ${\mathbf{ i}} = \zeta(a+a')$;  we want to show
${\mathbf{ j}}={\mathbf{ i}}$.  From Lemma \ref{kernel.biperm.eca.1}{\bf(c)}, it suffices
to show that $i_0=j_0$.  But the operation on ${\mathcal{ K}}$ is
componentwise addition.  Thus, $j_0 = k_0 + k'_0 = a + a' = i_0$.
Hence, $\zeta$ is a homomorphism; \ being
bijective, $\zeta$ is thus an isomorphism.  All other claims follow.
 {\tt \hrulefill $\Box$ } \end{list}  \medskip  

  Let $\eta$ be as in \S\ref{S:degree}, and
 for any ${\mathbf{ k}}\in{\mathcal{ K}}$, let ${\mathbf{ E}}_{\mathbf{ k}} \ = \ {\left\{ {\mathbf{ x}}\in{\mathcal{ A}}^{\mathbb{Z}} \; ; \; \eta({\mathbf{ x}}\bullet{\mathbf{ k}})>0 \right\} }$.  

\begin{lemma}{\sf \label{biperm.eca.3}}  
   \setcounter{enumi}{\thethm} \begin{list}{{\bf (\alph{enumii})}}{\usecounter{enumii}} 			{\setlength{\leftmargin}{0em} 			\setlength{\rightmargin}{0em}}
   \item $ {{{\boldsymbol{\sigma}}}^{}} ({\mathbf{ E}}_{\mathbf{ k}}) \ \displaystyle\raisebox{-0.6ex}{$\overline{\overline{{\scriptstyle{\mathrm{\mu}}}}}$} \ {\mathbf{ E}}_{ {{{\boldsymbol{\sigma}}}^{}} ({\mathbf{ k}})}$.

   \item Thus, if $ {{{\boldsymbol{\sigma}}}^{P}} ({\mathbf{ k}})={\mathbf{ k}}$, then 
$ {{{\boldsymbol{\sigma}}}^{P}} ({\mathbf{ E}}_{\mathbf{ k}})  \ \displaystyle\raisebox{-0.6ex}{$\overline{\overline{{\scriptstyle{\mathrm{\mu}}}}}$} \ {\mathbf{ E}}_{\mathbf{ k}}$.
\end{list}
 \end{lemma}
\bprf To prove {\bf(a)} it suffices to show
that $ {{{\boldsymbol{\sigma}}}^{}} ({\mathbf{ E}}_{\mathbf{ k}}) \ \raisebox{-1ex}{$\stackrel{\displaystyle\subset}{\scriptstyle{\mathrm{\mu}}}$} \ {\mathbf{ E}}_{ {{{\boldsymbol{\sigma}}}^{}} ({\mathbf{ k}})}$ (and
then, by symmetric reasoning, that ${\mathbf{ E}}_{ {{{\boldsymbol{\sigma}}}^{}} ({\mathbf{ k}})} \ \raisebox{-1ex}{$\stackrel{\displaystyle\subset}{\scriptstyle{\mathrm{\mu}}}$} \
 {{{\boldsymbol{\sigma}}}^{}} ({\mathbf{ E}}_{\mathbf{ k}})$.)  To show this, we define the measure
$\mu_{({\mathbf{ k}})}$ by $\mu_{({\mathbf{ k}})}({\mathbf{ U}}) \ = \
\mu\left({\mathbf{ E}}_{\mathbf{ k}}\cap( {\mathbf{ U}} \bullet{\mathbf{ k}}^{-1}) \right)$.  Then $\mu_{({\mathbf{ k}})}$
is absolutely continuous with respect to $\mu$ (by reasoning
similar to Lemma \ref{mu.n.ll.mu}).
 Lemma \ref{eta.shift.inv}{\bf(a)} says $\eta$ is $ {{{\boldsymbol{\sigma}}}^{}} $-invariant {\rm ($\mu$-\ae)};
hence  $\eta$ is $ {{{\boldsymbol{\sigma}}}^{}} $-invariant ($\mu_{({\mathbf{ k}})}$-\ae),
by reasoning similar to Corollary \ref{eta.constant.2}.
Thus, for ${{\forall}_{\mu}\;} {\mathbf{ x}}\in{\mathbf{ E}}_{\mathbf{ k}}$, \quad 
$0 \ < \  \eta({\mathbf{ x}}\bullet{\mathbf{ k}}) \  =  \   \eta\left(\rule[-0.5em]{0em}{1em} {{{\boldsymbol{\sigma}}}^{}} ({\mathbf{ x}}\bullet{\mathbf{ k}})\right)
\  =  \   \eta\left(\rule[-0.5em]{0em}{1em} {{{\boldsymbol{\sigma}}}^{}} ({\mathbf{ x}}) \bullet  {{{\boldsymbol{\sigma}}}^{}} ({\mathbf{ k}})\right)$, \
and thus, $ {{{\boldsymbol{\sigma}}}^{}} ({\mathbf{ x}})\in{\mathbf{ E}}_{ {{{\boldsymbol{\sigma}}}^{}} ({\mathbf{ k}})}$.  Hence  $ {{{\boldsymbol{\sigma}}}^{}} ({\mathbf{ E}}_{\mathbf{ k}}) \ \raisebox{-1ex}{$\stackrel{\displaystyle\subset}{\scriptstyle{\mathrm{\mu}}}$} \ {\mathbf{ E}}_{ {{{\boldsymbol{\sigma}}}^{}} ({\mathbf{ k}})}$. {\tt \hrulefill $\Box$ } \end{list}  \medskip  

  Recall from  \S\ref{S:degree} that ${\mathcal{ E}}({\mathbf{ x}}) \ = \ {\left\{ {\mathbf{ y}}\in{\mathcal{ F}}\left({\mathbf{ x}}\right) \; ; \; \eta({\mathbf{ y}})>0 \right\} }$.

\begin{cor}{\sf \label{biperm.eca.subkernel}}  
 If $\mu$ is $\Phi$-invariant and
totally $ {{{\boldsymbol{\sigma}}}^{}} $-ergodic, then there is a $ {{{\boldsymbol{\sigma}}}^{}} $-invariant subgroup
${\mathcal{ J}}\subset{\mathcal{ K}}$ so that, for ${{\forall}_{\mu}\;} {\mathbf{ x}}\in{\mathcal{ A}}^{\mathbb{Z}}$, \quad
${\mathcal{ E}}({\mathbf{ x}}) \ = \ {\mathbf{ x}}\bullet{\mathcal{ J}}$.
 \end{cor}
\bprf  Define ${\mathcal{ J}}={\left\{ {\mathbf{ k}} \in{\mathcal{ K}} \; ; \; \mu({\mathbf{ E}}_{\mathbf{ k}})>0 \right\} }$.  

\refstepcounter{claimcount}                {\bf Claim \theclaimcount: \ }{\sl  For any ${\mathbf{ j}}\in{\mathcal{ J}}$,\quad $\mu({\mathbf{ E}}_{\mathbf{ j}})=1$.}
\bprf
  By Lemma \ref{kernel.biperm.eca.1}{\bf(e)}, find $P\in{\mathbb{N}}$ so that
$ {{{\boldsymbol{\sigma}}}^{P}} ({\mathbf{ j}})={\mathbf{ j}}$.  But then Lemma \ref{biperm.eca.3}{\bf(b)}
says that $ {{{\boldsymbol{\sigma}}}^{P}} ({\mathbf{ E}}_{\mathbf{ j}})={\mathbf{ E}}_{\mathbf{ j}}$.  But $\mu$ is
$ {{{\boldsymbol{\sigma}}}^{P}} $-ergodic, so this means that $\mu({\mathbf{ E}}_{\mathbf{ j}})=1$.
 {\tt \dotfill~$\Diamond$~[Claim~\theclaimcount] }\end{list}

\refstepcounter{claimcount}                {\bf Claim \theclaimcount: \ }{\sl  For ${{\forall}_{\mu}\;}{\mathbf{ x}}\in{\mathcal{ A}}^{\mathbb{Z}}$, \quad  ${\mathcal{ E}}({\mathbf{ x}})={\mathbf{ x}}\bullet{\mathcal{ J}}$.}
\bprf
First note that ${\mathcal{ E}}({\mathbf{ x}}) \ = \ {\left\{ {\mathbf{ x}}\bullet{\mathbf{ k}} \; ; \; {\mathbf{ k}}\in{\mathcal{ K}}, \ {\mathbf{ x}}\in{\mathbf{ E}}_{\mathbf{ k}} \right\} }$.  Thus, we want to show that, for ${{\forall}_{\mu}\;}{\mathbf{ x}}\in{\mathcal{ A}}^{\mathbb{Z}}$, and
all ${\mathbf{ k}}\in{\mathcal{ K}}$,\quad $\left( \ \rule[-0.5em]{0em}{1em}       \begin{minipage}{40em}       \begin{tabbing}         ${\mathbf{ x}}\in{\mathbf{ E}}_{\mathbf{ k}}$        \end{tabbing}      \end{minipage} \ \right)
\iff\left( \ \rule[-0.5em]{0em}{1em}       \begin{minipage}{40em}       \begin{tabbing}         ${\mathbf{ k}}\in{\mathcal{ J}}$        \end{tabbing}      \end{minipage} \ \right)$.  Observe that
 $\mu\left(\displaystyle\bigcup_{{\mathbf{ k}}\in{\mathcal{ K}}\setminus{\mathcal{ J}}} {\mathbf{ E}}_{\mathbf{ k}}\right) \ = \ 0$ \
(by definition of ${\mathcal{ J}}$) \ 
and \  $\mu\left(\displaystyle\bigcap_{{\mathbf{ j}}\in{\mathcal{ J}}} {\mathbf{ E}}_{\mathbf{ j}}\right) \ = \ 1$ \
(by Claim 1).

Thus, ${\mathbf{ x}} \in \displaystyle \bigcap_{{\mathbf{ j}}\in{\mathcal{ J}}} {\mathbf{ E}}_{\mathbf{ j}}  \ \setminus \  \bigcup_{{\mathbf{ k}}\in{\mathcal{ K}}\setminus{\mathcal{ J}}} {\mathbf{ E}}_{\mathbf{ k}}$ for  ${{\forall}_{\mu}\;}{\mathbf{ x}}\in{\mathcal{ A}}^{\mathbb{Z}}$.  The claim follows.
 {\tt \dotfill~$\Diamond$~[Claim~\theclaimcount] }\end{list}

  Let ${\mathcal{ U}}= {\mathbf{ E}}_{\mathbf{ e}}={\left\{ {\mathbf{ x}}\in{\mathcal{ A}}^{\mathbb{Z}} \; ; \; \eta({\mathbf{ x}})>0 \right\} }$.

\refstepcounter{claimcount}                {\bf Claim \theclaimcount: \ }{\sl  If \ ${\mathbf{ k}}\in{\mathcal{ K}}$, then
\  $\left( \ \rule[-0.5em]{0em}{1em}       \begin{minipage}{40em}       \begin{tabbing}         ${\mathbf{ k}}\in{\mathcal{ J}}$        \end{tabbing}      \end{minipage} \ \right)\iff
\left( \ \rule[-0.5em]{0em}{1em}       \begin{minipage}{40em}       \begin{tabbing}         ${\mathcal{ U}}\bullet{\mathbf{ k}} \ \subset \ {\mathcal{ U}}$,  modulo a set of measure zero        \end{tabbing}      \end{minipage} \ \right)$.}
\bprf
  Claim 1 implies that $\mu({\mathcal{ U}})=1$.  Thus,
\begin{eqnarray*}
\left( \ \rule[-0.5em]{0em}{1em}       \begin{minipage}{40em}       \begin{tabbing}         ${\mathbf{ k}}\in{\mathcal{ J}}$        \end{tabbing}      \end{minipage} \ \right)&\Leftarrow\!\mbox{\tiny{$\mathrm{C1}$}}\!\Rightarrow&
\left( \ \rule[-0.5em]{0em}{1em}       \begin{minipage}{40em}       \begin{tabbing}         $\mu\left({\mathcal{ U}}\cap{\mathbf{ E}}_{\mathbf{ k}}\right)\ = \ 1$        \end{tabbing}      \end{minipage} \ \right)
\quad\Leftarrow\!\mbox{\tiny{$\mathrm{DE}$}}\!\Rightarrow\quad
\left( \ \rule[-0.5em]{0em}{1em}       \begin{minipage}{40em}       \begin{tabbing}         For ${{\forall}_{\mu}\;} {\mathbf{ u}}\in{\mathcal{ U}}$, \quad $\eta({\mathbf{ u}}\bullet{\mathbf{ k}})>0$        \end{tabbing}      \end{minipage} \ \right)
\\&\Leftarrow\!\mbox{\tiny{$\mathrm{DU}$}}\!\Rightarrow&
\left( \ \rule[-0.5em]{0em}{1em}       \begin{minipage}{40em}       \begin{tabbing}         For ${{\forall}_{\mu}\;} {\mathbf{ u}}\in{\mathcal{ U}}$, \quad ${\mathbf{ u}}\bullet{\mathbf{ k}}\in{\mathcal{ U}}$        \end{tabbing}      \end{minipage} \ \right)
\quad\iff\quad
\left( \ \rule[-0.5em]{0em}{1em}       \begin{minipage}{40em}       \begin{tabbing}         ${\mathcal{ U}}\bullet{\mathbf{ k}} \ \raisebox{-1ex}{$\stackrel{\displaystyle\subset}{\scriptstyle{\mathrm{\mu}}}$} \ {\mathcal{ U}}$        \end{tabbing}      \end{minipage} \ \right).
\end{eqnarray*}
 {\bf(C1)} is Claim 1. {\bf(DE)} is by definition of ${\mathbf{ E}}_{\mathbf{ k}}$.
 {\bf(DU)} is by definition of ${\mathcal{ U}}$.
 {\tt \dotfill~$\Diamond$~[Claim~\theclaimcount] }\end{list}

\refstepcounter{claimcount}                {\bf Claim \theclaimcount: \ }{\sl  ${\mathcal{ J}}$ is a subgroup of ${\mathcal{ K}}$.}
\bprf
Let ${\mathbf{ j}}_1,{\mathbf{ j}}_2\in{\mathcal{ J}}$, and let ${\mathbf{ j}} = {\mathbf{ j}}_1\bullet {\mathbf{ j}}_2$.  Then
 Claim 3 says \
$ {\mathcal{ U}}\bullet {\mathbf{ j}}
\ = \
 ({\mathcal{ U}}\bullet{\mathbf{ j}}_1)\bullet{\mathbf{ j}}_2
\ \subset\
{\mathcal{ U}}\bullet{\mathbf{ j}}_1
\ \subset\ {\mathcal{ U}}$,
\quad (modulo sets of measure zero).
Thus, Claim 3 implies ${\mathbf{ j}}\in{\mathcal{ J}}$ also.  Hence, ${\mathcal{ J}}$ is closed under
`$\bullet$'.  Since ${\mathcal{ J}}$ is finite, it is a subgroup.
 {\tt \dotfill~$\Diamond$~[Claim~\theclaimcount] }\end{list}

It remains to show that $ {{{\boldsymbol{\sigma}}}^{-1}} ({\mathcal{ J}})={\mathcal{ J}}$.  To see this, 
let ${\mathbf{ k}}\in{\mathcal{ K}}$.  Then
\begin{eqnarray*}
\left( \ \rule[-0.5em]{0em}{1em}       \begin{minipage}{40em}       \begin{tabbing}         ${\mathbf{ k}}\in{\mathcal{ J}}$        \end{tabbing}      \end{minipage} \ \right)
&\Leftarrow\!\mbox{\tiny{$\mathrm{C1}$}}\!\Rightarrow&\left( \ \rule[-0.5em]{0em}{1em}       \begin{minipage}{40em}       \begin{tabbing}         $\mu({\mathbf{ E}}_{\mathbf{ k}})=1$        \end{tabbing}      \end{minipage} \ \right)
\quad\Leftarrow\!\mbox{\tiny{$\mathrm{(*)}$}}\!\Rightarrow\quad
 \left( \ \rule[-0.5em]{0em}{1em}       \begin{minipage}{40em}       \begin{tabbing}         $\mu\left(\rule[-0.5em]{0em}{1em} {{{\boldsymbol{\sigma}}}^{-1}} ({\mathbf{ E}}_{\mathbf{ k}})\right)=1$        \end{tabbing}      \end{minipage} \ \right)\\
&\Leftarrow\!\mbox{\tiny{$\mathrm{(\dagger)}$}}\!\Rightarrow& \left( \ \rule[-0.5em]{0em}{1em}       \begin{minipage}{40em}       \begin{tabbing}         $\mu\left({\mathbf{ E}}_{ {{{\boldsymbol{\sigma}}}^{}} ({\mathbf{ k}})}\right)=1$        \end{tabbing}      \end{minipage} \ \right)
\ \ \Leftarrow\!\mbox{\tiny{$\mathrm{C1}$}}\!\Rightarrow\ \  \left( \ \rule[-0.5em]{0em}{1em}       \begin{minipage}{40em}       \begin{tabbing}         $ {{{\boldsymbol{\sigma}}}^{}} ({\mathbf{ k}})\in{\mathcal{ J}}$        \end{tabbing}      \end{minipage} \ \right)
\ \ \iff\ \
 \left( \ \rule[-0.5em]{0em}{1em}       \begin{minipage}{40em}       \begin{tabbing}         ${\mathbf{ k}}\in {{{\boldsymbol{\sigma}}}^{-1}} ({\mathcal{ J}})$        \end{tabbing}      \end{minipage} \ \right).
\end{eqnarray*}
 Here, {\bf(C1)} is by Claim 1, 
 $(*)$ is because $\mu$ is $ {{{\boldsymbol{\sigma}}}^{}} $-invariant, and
$(\dagger)$ is by Lemma \ref{biperm.eca.3}{\bf(a)}.
 {\tt \hrulefill $\Box$ } \end{list}  \medskip  

\begin{cor}{\sf \label{biperm.eca.subkernel.2}}  
 Let $J=|{\mathcal{ J}}|$.  Then  $\entr{\Phi} = \log(J)$, and $\Phi$ is $J$-to-1
{\rm ($\mu$-\ae)}.
 \end{cor}
\bprf
  Combine Corollary \ref{biperm.eca.subkernel} with 
Corollary \ref{eta.constant.5}.
 {\tt \hrulefill $\Box$ } \end{list}  \medskip  

\bprf[Proof of Theorem \ref{biperm.eca.measure}]
  If $\entr{\Phi}>0$, then Corollary \ref{biperm.eca.subkernel.2} says
$|{\mathcal{ J}}|>1$, so ${\mathcal{ J}}$ is a nontrivial $ {{{\boldsymbol{\sigma}}}^{}} $-invariant subgroup of
${\mathcal{ K}}$.  Thus, ${\mathcal{ J}}={\mathcal{ K}}$, which means $|{\mathcal{ J}}| \ = \ |{\mathcal{ K}}| \ = \ |{\mathcal{ A}}|$,
where the second equality is by Lemma \ref{kernel.biperm.eca.1}{\bf(d)}.
Thus, $\entr{\Phi}=\log|{\mathcal{ A}}|$.  Thus, $\entr{ {{{\boldsymbol{\sigma}}}^{}} }=\log|{\mathcal{ A}}|$,
which means $\mu$ must be the uniform measure.
 {\tt \hrulefill $\Box$ } \end{list}  \medskip  

  Lemmas \ref{kernel.biperm.eca.1}{\bf(g)}  and 
\ref{kernel.biperm.eca.2}{\bf(e)}  provide conditions under which ${\mathcal{ K}}$
has no $ {{{\boldsymbol{\sigma}}}^{}} $-invariant subgroups.  
For example, suppose $p\in{\mathbb{N}}$ is prime, and let ${\mathcal{ A}} = ({{\mathbb{Z}}_{/p}})^N$
for some $N>0$.  Then ${\mathcal{ A}}$ is a vector space over the field ${{\mathbb{Z}}_{/p}}$,
and $\rho:{\mathcal{ A}}{{\longrightarrow}}{\mathcal{ A}}$ is a group automorphism iff $\rho$ is a
${{\mathbb{Z}}_{/p}}$-linear automorphism.  Thus, $\rho$ can be described by
an $N\times N$ matrix ${\mathbf{ M}}$ of coefficients in ${{\mathbb{Z}}_{/p}}$.  Furthermore,
${\mathcal{ B}}\subset{\mathcal{ A}}$ is a ($\rho$-invariant) subgroup iff ${\mathcal{ B}}$ is
a ($\rho$-invariant) subspace.  The structure of
$\rho$-invariant subspaces in ${\mathcal{ A}}$ is described by the {\em rational
canonical form} of $\rho$;  this is an $N \times N$ matrix ${\widetilde{\mathbf{ M}}}$,
similar to ${\mathbf{ M}}$, having the block-diagonal form

\[
  {\widetilde{\mathbf{ M}}} \quad=\quad 
\left[\begin{array}{ccc}
{\mathbf{ M}}_1 & \ldots & 0 \\
 \vdots &\ddots &\vdots \\
    0 &  \ldots & {\mathbf{ M}}_L
\end{array}
\right]
\]
where each {\bf component matrix} ${\mathbf{ M}}_\ell$ has the form:
$ {\mathbf{ M}}_\ell \ = \ $
{\small $
\left[\begin{array}{cccccc}
0 & 0 & 0 &\ldots & 0 & m_1\\
 1 & 0 & 0 & \ldots & 0 & m_2 \\
 0 & 1 & 0 & \ldots & 0 & m_3 \\
\vdots &\vdots & \vdots &\ddots &\vdots  &\vdots  \\
 0 & 0 & 0 & \ldots & 0 & m_{r-1} \\
 0 & 0 & 0 & \ldots & 1& m_{r} 
\end{array}
\right]$,}
 for some $r>0$ and $m_1,\ldots,m_r\in{{\mathbb{Z}}_{/p}}$.
  Each component matrix corresponds to a $\rho$-invariant subspace of
${\mathcal{ A}}$.  If ${\widetilde{\mathbf{ M}}}$ has only one component, then we say ${\widetilde{\mathbf{ M}}}$ is {\bf
simple.}  We call the automorphism $\rho$ {\bf simple} if its rational
canonical form is simple.  It follows:

\begin{lemma}{\sf \label{simple.rat.can.form}}  
$\left( \ \rule[-0.5em]{0em}{1em}       \begin{minipage}{40em}       \begin{tabbing}          $\rho$ is simple        \end{tabbing}      \end{minipage} \ \right)
\iff
\left( \ \rule[-0.5em]{0em}{1em}       \begin{minipage}{40em}       \begin{tabbing}          ${\mathcal{ A}}$ has no nontrival $\rho$-invariant subspaces.        \end{tabbing}      \end{minipage} \ \right)$.
\hrulefill\ensuremath{\Box}
 \end{lemma}

\begin{cor}{\sf }  
  Let ${\mathcal{ A}} = ({{\mathbb{Z}}_{/p}})^N$ and let ${\mathcal{ A}}^{\mathbb{Z}}$ 
be the product group. Let $\Phi:{\mathcal{ A}}^{\mathbb{Z}}{{\longrightarrow}}{\mathcal{ A}}^{\mathbb{Z}}$ be
a bipermutative ECA with local map $\phi(a_0,a_1) =
\phi_0(a_0) + \phi_1(a_1)$.  Suppose $\rho=-\phi_1^{-1}\circ \phi_0$
is simple.  Then the conclusion of Theorem \ref{biperm.eca.measure}
holds.
 \end{cor}
\bprf
Combine Lemma \ref{simple.rat.can.form} 
with parts {\bf(c)} and {\bf(e)} 
of  Lemma \ref{kernel.biperm.eca.2}.
 {\tt \hrulefill $\Box$ } \end{list}  \medskip  

        \refstepcounter{thm}                     \begin{list}{} 			{\setlength{\leftmargin}{1em} 			\setlength{\rightmargin}{0em}}        \item {\bf Example \thethm:} 
 Let ${\mathcal{ A}} = ({{\mathbb{Z}}_{/7}})^4$, and suppose $\phi(a_0,a_1) = \phi_0(a_0) + 
a_1$, where $\phi_0$ has matrix 
{\small \[
\left[\begin{array}{cccc}
0 & 0 & 0 & 1 \\
1 & 0 & 0 & 1 \\
0 & 1 & 0 & 1 \\
0 & 0 & 1 & 1 \\
\end{array}
\right]
\] }
Thus, $\rho=-\phi_0$ is simple.  Hence, if $\mu$ is $\Phi$-invariant
and totally $ {{{\boldsymbol{\sigma}}}^{}} $-ergodic, and $\entr{\Phi}>0$, then $\mu$ is the
uniform measure.                   \hrulefill  \end{list}   			

\subsection*{Conclusion}

  We have characterized the invariant measures for several natural
families of bipermutative cellular automata.  Many questions remain
unanswered.  For example, in \S\ref{S:group} and \S\ref{S:eca}, we
exploited an algebraic structure on ${\mathcal{ A}}^{\mathbb{Z}}$ to study the
$\Phi$-invariant measures.  What other algebraic properties of the quasigroup
structure of ${\mathcal{ A}}$ can be exploited in this way?
 
  Also, if $\Phi$ is a QGCA on
${\mathcal{ A}}^{\mathbb{Z}}$, then the system $({\mathcal{ A}}^{\mathbb{Z}},\Phi,\lambda)$ is measurably
isomorphic to the uniform Bernoulli shift $({\mathcal{ A}}^{\mathbb{N}}, {{{\boldsymbol{\sigma}}}^{}} ,\lambda)$
\cite{Shereshevsky}.  Theorem \ref{loop.ca.K.2.1} suggests
that, if $\mu$ is any positive-entropy, $\Phi$-invariant measure,
then  the system $({\mathcal{ A}}^{\mathbb{Z}},\Phi,\mu)$ is 
isomorphic to $({\mathcal{ K}}^{\mathbb{Z}}, {{{\boldsymbol{\sigma}}}^{}} ,\kappa)$,
where ${\mathcal{ K}}$ is an alphabet of $K$ letters and $\kappa$ is the
uniform Bernoulli measure on ${\mathcal{ K}}^{\mathbb{N}}$.  Is this true?

 Finally, Example \ref{counterexample} refuted Conjecture
\ref{subloop.conjecture}, but did so by using a structural
decomposition ${\mathcal{ A}}={\mathcal{ C}}\times{\mathcal{ Q}}$ to get an invariant measure without
full support.  This leaves us with the following:

\paragraph*{Conjecture:}
{\em Let  $({\mathcal{ A}},*)$ be a quasigroup and let
 $\Phi:{\mathcal{ A}}^{\mathbb{Z}}{{\longrightarrow}}{\mathcal{ A}}^{\mathbb{Z}}$ be the corresponding QGCA.  Let
$\mu$ be a $\Phi$-invariant and $ {{{\boldsymbol{\sigma}}}^{}} $-ergodic measure.  

 1.\ If $\mu$ has full support, then $\mu=\lambda$. 

 2.\ If $({\mathcal{ A}},*)$ is simple  (ie. has no nontrivial quotients), and
$\entr{\Phi}>0$, then $\mu=\lambda$.
}

{\footnotesize
\bibliographystyle{plain}
\bibliography{bibliography}
}
\end{document}